\documentclass[12pt]{article}
\usepackage{amsmath,amsfonts,amsthm,amssymb,eucal,amsxtra,amscd}
\newcommand{\tp}{\widehat{\otimes}}
\newcommand{\K}{\mathbb{K}}
\newcommand{\D}{\Hat \Delta}
\newcommand{\A}{\Hat A}
\newcommand{\AAA}{\Hat{\Hat A}}

\newcommand{\id}{\operatorname{id}}
\newcommand{\RR}{\mathbf R}

\begin{document}

\newtheorem{teo}{Theorem}[section]
\newtheorem{lem}[teo]{Lemma}
\newtheorem{prop}[teo]{Proposition}
\newtheorem{cor}[teo]{Corollary}

\numberwithin{equation}{section}
\pagestyle{plain}
\title{Non-Archimedean Duality: Algebras, Groups, and Multipliers}
\author{Anatoly N. Kochubei\\
\footnotesize Institute of Mathematics,\\
\footnotesize National Academy of Sciences of Ukraine,\\
\footnotesize Tereshchenkivska 3, Kiev, 01601 Ukraine
\\ \footnotesize E-mail: \ kochubei@i.com.ua}
\date{}
\maketitle

\bigskip
\begin{abstract}
We develop a duality theory for multiplier Banach-Hopf algebras over a non-Archime\-dean field $\K$. As examples, we consider algebras corresponding to discrete groups and zero-dimensional locally compact groups with $\K$-valued Haar measure, as well as algebras of operators generated by regular representations of discrete groups.
\end{abstract}


\bigskip
\section{INTRODUCTION}

Non-Archimedean harmonic analysis, the Fourier analysis of functions $f:\ G\to \K$ where $G$ is a group, $\K$ is a non-Archimedean valued field, was initiated in the thesis \cite{Sch1} by W. H. Schikhof; this subject should not be confused with the study of complex-valued functions on non-Archimedean structures started in another great thesis, by J. Tate. For Abelian groups admitting a $\K$-valued Haar measure (this class was described by Monna and Springer \cite{MS}), Schikhof proved an analog of Pontryagin's duality theorem; see also \cite{vR,vRS}. Duality theorems for compact groups and related group algebras in the non-Archimedean setting were proved later by Schikhof \cite{Sch2} and Diarra \cite{D74}; for related subjects see \cite{D96}.

In this paper we develop a duality theory for general multiplier Banach-Hopf algebras over $\K$, including a description of a dual object and a biduality theorem. This covers three main examples: 1) the algebra $A_1=c_0 (G_1)$ of $\K$-valued functions on a discrete group $G_1$ tending to zero by the filter of complements of finite sets, with pointwise operations; 2) the algebra $A_2=C_0(G_2)$ of $\K$-valued continuous functions on a zero-dimensional locally compact group $G_2$ admitting a $\K$-valued Haar measure, and 3) the algebra $A_3=A(G_3)$ generated by the right regular representation of a discrete group $G_3$. The first two examples describe, as much as possible, the group duality (groups without a $\K$-valued Haar measure possess quite pathological properties \cite{D79} and are not accessible so far).

The third algebra $A(G_3)$ appeared \cite{K14} as a group algebra belonging, for some classes of groups, to the class of operator algebras with Baer reductions, a class seen as a non-Archimedean counterpart of von Neumann algebras \cite{K13}. By its properties \cite{K14}, this algebra can be seen simultaneously as an analog of the reduced group $C^*$-algebra. In this paper we see its third face -- as a Banach-Hopf algebra, it is isomorphic to the object dual to $A_1$ and can be considered as a non-Archimedean compact quantum group. For another approach to non-Archimedean quantum groups see \cite{So}.

Our methods and framework are very close to the purely algebraic setting developed by Van Daele \cite{VD94,VD98}; in the group situation this theory deals with complex-valued functions on a group different from zero at finitely many points. It is interesting that the transition from the algebraic duality theory to the non-Archimedean topological theory is easier than the similar, in principle, transitions to the cases of $C^*$-algebras and von Neumann algebras over $\mathbb C$ \cite{KV00,KV03,VVD,VD14}. The main reason is the fact that, for example, $c_0 (G_1)$ is simultaneously the basic Banach space; no Hilbert space is needed, and a reasonable analog of the latter does not exist in the non-Archimedean case where there are no natural involutions. In contrast to the classical situations, in all the above examples (and in our general results) the counit and antipode are everywhere defined bounded mappings.

The contents of the paper are as follows. In Section 2, we present some preliminaries from non-Archimedean analysis; we also discuss the notion of a multiplier and properties of multiplier algebras. In Section 3, we define and study non-Archimedean multiplier Banach-Hopf algebras. Section 4 is devoted to the dual object and the biduality theorem. In Section 5, we describe in detail the above three examples.

The proofs of many results about non-Archimedean multiplier Banach-Hopf algebras are very similar to their classical counterparts; in such cases we will just refer refer to papers by Van Daele \cite{VD94,VD98} and their expositions in \cite{Ti,Sun}.

\medskip
{\bf Acknowledgements.} The author is grateful to Leonid Vainerman for very helpful consultations, and to the anonymous referee for valuable remarks. I am also pleased to acknowledge the influence of Georgiy Kac (1924-1978) whose seminar I attended in 1970s, and that of Wim Schikhof (1937-2014).

\section{Preliminaries}

{\bf 2.1.} {\it Non-Archimedean Banach spaces} \cite{BGR,PGS,vR}. A non-Archimedean valued field is a field $\K$ with a nonnegative real absolute value (or valuation) $|\lambda |$, $\lambda \in \K$, such that $|\lambda |=0$ if and only if $\lambda =0$, $|\lambda \mu|=|\lambda|\cdot |\mu|$, $|\lambda +\mu|\le \max (|\lambda|, |\mu|)$, $\lambda, \mu \in \K$. Below we assume that the valuation is nontrivial, that is we exclude the case where $|\lambda |=1$ for every $\lambda \ne 0$. The field with a nontrivial valuation is a nondiscrete totally disconnected topological field with respect to the topology induced by the ultrametric $(\lambda ,\mu )\mapsto |\lambda -\mu|$. Typical examples are the field $\mathbb Q_p$ of $p$-adic numbers ($p$ is a prime number) and the field $\mathbb C_p$, the completion of an algebraic closure of $\mathbb Q_p$ with respect to a valuation obtained by extension of the one from $\mathbb Q_p$. On $\mathbb Q_p$, the valuation is discrete -- the absolute value takes values from the set $\{ p^n,n\in \mathbb Z\}$, while the valuation on $\mathbb C_p$ is dense, its possible values are $p^\nu$, $\nu \in \mathbb Q$. Many properties of these fields are quite different.

Let $E$ be a $\K$-vector space. A {\it norm} on $E$ is a map $\| \cdot \|:\ E\to [0,\infty )$, such that $\| x\|=0$ if and only if $x=0$; $\| \lambda x\|=|\lambda|\cdot \| x\|$; $\| x+y\| \le \max (\| x\|,\|y\|)$ for all $x,y\in E$, $\lambda \in \K$. $E$ is called a (non-Archimedean) Banach space, if $E$ is complete with respect to the ultrametric $(x,y)\mapsto \|x-y\|$. Under this ultrametric, all the sets $\{ x\in E:\ \|x-x_0\| \le c\}$, $\{ x\in E:\ \|x-x_0\| <c\}$, $\{ x\in E:\ \|x-x_0\| =c\}$ ($x_0\in e,c>0$) are both open and closed. A normed space $E$ is a Banach space, if and only if every sequence converging to zero is summable.

The above phenomena are of purely non-Archimedean nature. Meanwhile most of the notions and results of classical functional analysis have their non-Archimedean counterparts. Here we will touch the notions of non-Archimedean orthogonality, separation property and tensor product.

We say that vectors $x$ and $y$ of a normed space $E$ are {\it orthogonal}, if, for all $\lambda,\mu \in \K$,
\begin{equation}
\| \lambda x+\mu y\| =\max (\| \lambda x\|,\|\mu y\|).
\end{equation}
If $\|x\|=\|y\|=1$, the equality (2.1) turns into the orthonormality relation
$$
\| \lambda x+\mu y\| =\max (| \lambda |,|\mu |).
$$
In an obvious way, this property is extended to any finite number of elements.

There is also a weaker property of $t$-orthogonality of a system $\{ e_1,\ldots ,e_n\}$, $0<t\le 1$, defined by the property that
$$
\| \lambda_1e_1+\cdots +\lambda_ne_n\|\ge t\max\limits_{1\le i\le n}\|\lambda_ie_i\|
$$
for any $\lambda_1,\ldots ,\lambda_n\in \K$. The 1-orthogonality coincides with the above orthogonality.

A sequence $e_1,e_2,\ldots \in E$ is called a $t$-orthogonal basis, if any its finite subset is $t$-orthogonal and every $x\in E$ has a convergent expansion $x=\sum\limits_{n=1}^\infty \lambda_ne_n$, $\lambda_n\in \K$. The basis is called {\it orthogonal}, if $t=1$, and {\it orthonormal} if, in addition, $\|e_n\|=1$ for all $n\in \mathbb N$. See \cite{BGR,PGS,vR} regarding classes of Banach spaces possessing orthonormal bases.

It was a nontrivial problem whether an analog of the Hahn-Banach theorem holds in the non-Archimedean case. We will not need its solution in this paper (see \cite{PGS,vR}); note only that it holds for spaces over $\mathbb Q_p$ and typically does not hold for spaces over $\mathbb C_p$. However we will need the following weaker property resembling some classical applications of the Hahn-Banach theorem.

A Banach space $E$ over $\K$ is called {\it polar}, if it has the following property. Let $B$ be the unit ball in $E$, that is $B=\{ x\in E:\ \|x\|\le 1\}$, and let $y\notin B$. Then there exists such a linear continuous functional $f:\ E\to \K$ that $|f(x)|\le 1$ for $x\in B$ and $|f(y)|>1$; see Lemma 4.4.4 in \cite{PGS} regarding the equivalence of this definition to some others.

If a Banach space $E$ is polar, then linear continuous functionals separate points of $E$. All the spaces encountered in this paper are polar (see Section 4.4 in \cite{PGS}).

Let $E$ and $F$ be Banach spaces over $\K$. On the algebraic tensor product $E\otimes F$, we define a norm setting
\begin{equation}
\| g\| =\inf \left( \max\limits_{1\le i\le r}\|x_i\| \cdot \|y_i\|\right)
\end{equation}
where the infimum runs over all possible representations
$$
g=\sum\limits_{i=1}^r x_i\otimes y_i,\quad x_i\in E,y_i\in F.
$$
The proof that (2.2) is a norm indeed, employs essentially the non-Archimedean properties of $E,F$; see \cite{PGS}, Corollary 10.2.10. Note also that $\|x\otimes y\|=\| x\|\cdot \|y\|$ (see \cite{vR}, Theorem 4.27).

The completion $E\tp F$ of $E\otimes F$ with respect to this norm is called {\it the completed tensor product} of $E$ and $F$.

\medskip
\begin{lem}
Let $w\in E\tp F$, $0<t<1$. Then there exists a $t$-orthogonal sequence $\{ a_i\}_{i=1}^\infty \subset E$, and a sequence $\{ b_i\}_{i=1}^\infty \subset F$, such that $\lim\limits_{i\to \infty}\| a_i\|\cdot \|b_i\|=0$,
\begin{equation}
w=\sum\limits_{i=1}^\infty a_i\otimes b_i
\end{equation}
and $\| w\| =\inf \left\{ \max\limits_i\| a_i\|\cdot \|b_i\|\right\}$ where the infimum runs over all possible representations (2.3). Once the $t$-orthogonal sequence $\{ a_i\}$ is chosen, the sequence $\{ b_i\}$ is unique and $t\sup\limits_{i\ge 1}\|a_i\| \cdot \|b_i\| \le \|w\| \le \sup\limits_{i\ge 1}\|a_i\| \cdot \|b_i\|$.
\end{lem}

\medskip
For the proof see \cite{vR}, Lemma 4.30 and Corollary 4.31, and \cite{PGS}, Corollary 10.2.10.

As in the classical case, given linear bounded mappings of Banach spaces, it is possible to define their tensor product acting on the completed tensor product of the spaces; see \cite{PGS}, Theorem 10.3.7.

\bigskip
{\bf 2.2.} {\it Banach algebras}. A $\K$-algebra $A$ is called a {\it Banach algebra}, if $A$ is simultaneously a Banach space with the norm $\|\cdot \|$, and $\|xy\|\le \|x\|\cdot \|y\|$, $x,y\in A$. If $A$ is an algebra with unit $e$, then we assume that $\|e\|=1$. We will often denote a unit by 1, if that does not lead to a confusion. We call a Banach algebra polar, if such is the underlying Banach space.

If $A$ and $B$ are Banach algebras, then $A\tp B$ is also a Banach algebra with the multiplication extending the relation $(a'\otimes b')(a''\otimes b'')=(a'a'')\otimes (b'b'')$; see \cite{BGR}, Section 3.1.1. If $A$ and $B$ are polar, then $A\tp B$ is polar too (\cite{PGS},Theorem 10.2.7).

We say that the product in an algebra $A$ is {\it nondegenerate} if the equality $ab=0$ for all $b$ implies $a=0$, and the equality $ab=0$ for all $a$ implies $b=0$.

\medskip
\begin{lem}
If $A$ and $B$ are polar Banach algebras with nondegenerate products, then the product in $A\tp B$ is nondegenerate.
\end{lem}

\medskip
{\it Proof}. Let $w\in A\tp B$, $(c\otimes d)w=0$ for all $c\in A,d\in B$. Writing $w$ in the form of (2.3), where the sequences $\{ a_i\},\{b_i\}$ are such as stated in Lemma 2.1, we find that
\begin{equation}
\sum\limits_{i=1}^\infty (ca_i)\otimes (db_i)=0.
\end{equation}

Let us apply to both sides of (2.4) the mapping $\id\otimes \chi$ where $\chi$ is a linear continuous functional on $B$. We get
$$
c\sum\limits_{i=1}^\infty a_i\chi (db_i)=0.
$$
Since $c$ is arbitrary and the product in $A$ is nondegenerate, we have $\sum\limits_{i=1}^\infty a_i\chi (db_i)=0$, and the $t$-orthogonality of $\{a_i\}$ implies the equalities $\chi (db_i)=0$ for all $i$. Since $\chi$ is arbitrary and $B$ is polar, we get $db_i=0$, thus $b_i=0$ and $w=0$. Similarly, if $w(c\otimes d)=0$ for any $c,d$, then $w=0$. $\qquad \blacksquare$

\medskip
We will often use the notion of a bounded approximate identity in a Banach algebra $A$ without a unit. By definition, a bounded approximate identity is a net $\{ e_\lambda\}_{\lambda \in \Lambda}\subset A$, contained within some ball, such that, for any $x\in A$, $e_\lambda x\to x$ and $xe_\lambda \to x$ in the topology of $A$.

\medskip
{\bf 2.3.} {\it Multipliers}. Let $A$ be a polar Banach algebra over $\K$ with nondegenerate product.

A {\it left multiplier} on $A$ is a linear continuous mapping $\rho_l:\ A\to A$, such that $\rho_l(ab)=\rho_l(a)b$ for all $a,b\in A$. Similarly, a {\it right multiplier} is defined by the property $\rho_r(ab)=a\rho_r(b)$, $a,b\in A$. A {\it multiplier} on $A$ is a couple $(\rho_l,\rho_r)$ of a left multiplier $\rho_l$ and a right multiplier $\rho_r$, such that
\begin{equation}
\rho_r(a)b=a\rho_l(b) \quad\text{for all $a,b\in A$}.
\end{equation}

If $\rho_l$ and $\rho_r$ are linear continuous mappings satisfying (2.5), then $(\rho_l,\rho_r)$ is a multiplier. Indeed, it follows from (2.5) that $\rho_r(ac)b=ac\rho_l(b)$ for any $a,b,c\in A$, and also that $\rho_r(c)b=c\rho_l(b)$, so that $a\rho_r(c)b=ac\rho_l(b)$. Subtracting we have
$$
[\rho_r(ac)-a\rho_r(c)]b=0 \quad \text{for all $b$}.
$$
It follows from the nondegeneracy that $\rho_r$ is a right multiplier. Similarly we show that $\rho_l$ is a left multiplier.

Denote by $L(A),R(A)$, and $M(A)$ the sets of left, right, and multipliers on $A$ respectively. The sets $L(A)$ and $R(A)$ are unital Banach algebras with respect to the usual norms of operators. $M(A)$ is a unital Banach algebra with the multiplication
$$
(\rho_l',\rho_r')(\rho_l'',\rho_r'')=(\rho_l'\rho_l'',\rho_r''\rho_r')
$$
and the norm $\|(\rho_l,\rho_r)\|=\max (\|\rho_l\|,\|\rho_r\|)$. Examples of multiplier algebras will be giver in Section 5 below.

For brevity, we will often write $ca=\rho_l(a)$ and $ac=\rho_r(a)$ for a multiplier $c=(\rho_l,\rho_r)$ and $a\in A$. In other words, in order to define a multiplier $c$, it suffices to define continuous multiplications $ac$ and $ca$ and to prove that $(ac)b=a(cb)$. This notation agrees with the obvious imbeddings $A\subset L(A)$, $A\subset R(A)$, and $A\subset M(A)$. When $A$ is unital, all the three multiplier algebras are isomorphic to $A$. Therefore below we consider non-unital algebras.

Let $A,B$ be Banach algebras over $\K$ with nondegenerate product. A continuous homomorphism $\Phi:\ A\to M(B)$ is called {\it nondegenerate} if the linear span of each of the sets $\{ \Phi (a)b:\ a\in A,b\in B\}$ and $\{ b\Phi (a):\ a\in A,b\in B\}$ is dense in $B$.

\medskip
\begin{teo}
If $A$ and $B$ have nondegenerate products, $A$ possesses a bounded approximate identity $\{e_\lambda\}_{\lambda \in \Lambda}$, and $\Phi:\ A\to M(B)$ is a nondegenerate continuous homomorphism, then $\Phi$ has a unique extension to a continuous homomorphism $\Phi_1:\ M(A)\to M(B)$.
\end{teo}

\medskip
{\it Proof}. Let $(L,R)\in M(A)$. Let us first define the operators $\Phi_1(L)$ and $\Phi_1(R)$ on linear sets $B_l=\operatorname{span}\{ \Phi (a)b:\ a\in A,b\in B\}$ and $B_r=\operatorname{span}\{ b\Phi (a):\ a\in A,b\in B\}$ respectively.

If $\Phi_1(L)$ is a required continuation onto $B_l$, then $\Phi (La)b=\Phi_1(L)\Phi (a)b$, $a\in A,b\in B$. This implies the uniqueness of the continuation. Similarly $b\Phi (aR)=b\Phi (a)\Phi_1(R)$, so that we get the uniqueness of the continuation $\Phi_1(R)$.

If there is a finite sum $\sum \Phi(a_i)b_i=0$, then for any $d\in A$, $e\in B$,
$$
e\Phi (d)\sum \Phi(La_i)b_i=e\sum \Phi(dLa_i)b_i=e\Phi (dR)\sum \Phi(a_i)b_i=0,
$$
so that $\sum \Phi(La_i)b_i=0$. Therefore it is legitimate to define $\Phi_1(L)$ setting
$$
\Phi_1(L)[\Phi (a)b]=\Phi (La)b.
$$
Similarly, $\Phi_1(R)$ is defined as follows:
$$
b\Phi (a)\Phi_1(R)=b\Phi (aR).
$$
Then we set $\Phi_1(L,R)=(\Phi_1(L),\Phi_1(R))$. Obviously, $\Phi_1$ is an extension of $\Phi$. Also,
$$
\Phi_1(L_1L_2)\Phi (a)b=\Phi (L_1L_2a)b=\Phi_1(L_1)\Phi (L_2a)b=\Phi_1(L_1)\Phi_1(L_2)\Phi (a)b.
$$
This relation, together with a similar one for right multipliers, shows that $\Phi_1$ is a homomorphism.

We have not yet checked that $\Phi_1$ sends multipliers to multipliers. Let
$$
z_1=\sum\limits_i\Phi (a_i)b_i,\quad z_2=\sum\limits_i\Phi (c_i)d_i.
$$
Then
\begin{multline*}
\Phi_1(L)(z_1z_2)=\Phi_1(L)\left( \sum\limits_{i,j}[\Phi (a_i)b_i]\cdot [\Phi (c_j)d_j]\right) =\sum\limits_{i,j}\Phi_1(L)(\Phi (a_i)(b_i(\Phi (c_j)d_j))\\
=\sum\limits_{i,j}\Phi (La_i)(b_i\cdot (\Phi (c_j)d_j))=\sum\limits_{i,j}[\Phi (La_i)b_i]\cdot \Phi (c_j)d_j=\Phi_1(L)(z_1)z_2.
\end{multline*}
Similarly we can get the required property for right multipliers.

Since $\Phi (Le_\lambda )\Phi (x)=\Phi (Le_\lambda x)\to \Phi (Lx)$ in $M(B)$ for every $x\in A$, we find that
$$
\Phi (Le_\lambda )z\longrightarrow \Phi_1(L)z \quad \text{for all $z\in B_l$}.
$$
The net $\{e_\lambda\}$ is bounded, so that $\|\Phi (Le_\lambda )z\|_B\le C\|z\|_B$ and $\|\Phi_1(L)z \|_B\le C\|z\|_B$. Thus, we have proved that $\Phi_1(L)$ is a bounded operator on $B_l$. Therefore it admits a unique extension to a bounded operator on $B$. Its multiplier property follows via the extension by continuity.

If $L_m\to 0$ in the uniform operator topology, then $L_me_\lambda \to 0$ uniformly in $\lambda \in \Lambda$. Therefore $\Phi (L_me_\lambda )\xi \to 0$ in $B$, uniformply with respect to $\lambda \in \Lambda$ and $\xi \in B$, $\|\xi\|\le 1$. We may take a limit in $\lambda$ and obtain that $\Phi_1 (L_m)\xi \to 0$ uniformly in $\xi \in B$, $\|\xi\|\le 1$, so that $\|\Phi_1 (L_m)\|\to 0$. This means the continuity of $\Phi_1$. The proof for right multipliers is similar. $\qquad \blacksquare$

\medskip
Similarly one can obtain the extension property for anti-homomorphisms, that is for mappings $\Phi$ satisfying the identily $\Phi (ab)=\Phi (b)\Phi (a)$, $a,b\in A$. The only difference is that an anti-homomorphism transforms a left multiplier into a right one, and vice versa:
$$
b\Phi (La)=b\Phi (a)\Phi_1(L)
$$
etc.

\medskip
{\it Remark}. The above proof follows \cite{VD94} and \cite{Jo}. A proof of the extension property given within a general theory of multipliers of complex Banach algebras \cite{Da} is based on the Cohen-Hewitt factorization theorem. It is not known whether the latter is valid in the non-Archimedean situation, but all its available proofs (see \cite{Ki}) fail in this case.

\bigskip
{\bf 2.4.} {\it Strict topology}. Below we will use systematically the multiplier algebra $M(A\tp A)$. The imbedding $A\subset M(A)$ implies the imbeddings $A\tp A\subset M(A)\tp M(A)\subset M(A\tp A)$.

Let us consider the subset $M_0(A\tp A)\subset M(A\tp A)$ consisting of such multipliers $x\in M(A\tp A)$ that $x(1\otimes a)$, $x(a\otimes 1)$, $(1\otimes a)x$, and $(a\otimes 1)x$ belong to $A\tp A$ for all $a\in A$. $M_0(A\tp A)$ is a $\K$-vector subspace of $M(A\tp A)$. The set of all seminorms $x\mapsto \|x(1\otimes a)\|_{A\tp A}$, $x\mapsto \|x(a\otimes 1)\|_{A\tp A}$, $x\mapsto \|(1\otimes a)x\|_{A\tp A}$, $x\mapsto \|(a\otimes 1)x\|_{A\tp A}$, $a\in A$ endow $M_0(A\tp A)$ with a locally convex topology, which will be called the {\it strict topology}.

\medskip
\begin{prop}
If $A$ possesses a bounded approximate identity $\{e_\lambda\}_{\lambda \in \Lambda}$, then $A\otimes A$ is dense in $M_0(A\tp A)$ in the strict topology. In particular, if $x\in M_0(A\tp A)$, then $x(1\otimes e_\lambda )\to x$, $x(e_\lambda \otimes 1)\to x$, $(e_\lambda \otimes 1)x\to x$, $(1\otimes e_\lambda )x\to x$ in the strict topology.
\end{prop}

\medskip
{\it Proof}. The first assertion follows from the second one, since, for example, $x(1\otimes e_\lambda)$ belongs to $A\tp A$ and can be approximated by elements from $A\otimes A$ in the topology of $A\tp A$, thus also in the strict topology.

If $y\in A\otimes A$, $y=\sum\limits_{i=1}^n p_i\otimes q_i$, $p_i,q_i\in A$, then $y-y(1\otimes e_\lambda )=\sum\limits_{i=1}^n p_i\otimes (q_i-q_ie_\lambda )$, so that
$$
\|y-y(1\otimes e_\lambda )\|_{A\tp A}\le \max\limits_i\left[ \|p_i\|\cdot \|q_i-q_ie\lambda \|\right]\to 0
$$
in $\lambda\in \Lambda$.

Let $\|e_\lambda\|\le C$. Then $\|(1\otimes e_\lambda\|_{M(A\tp A)}\le C$. If $X\in A\tp A$, then for any $\epsilon >0$, there exists $y\in A\otimes A$, such that$\|X-y\|_{A\tp A}<\min \left(\epsilon, \frac{\epsilon}C\right)$. We get
$$
\|X-X(1\otimes e_\lambda )\|_{A\tp A}\le \max\{ \epsilon,\|y-y(1\otimes e_\lambda )\|_{A\tp A},\|(X-y)(1\otimes e_\lambda )\|_{A\tp A}\}.
$$
Choosing $\lambda_0\in \Lambda$, such that $\|y-y(1\otimes e_\lambda )\|_{A\tp A}<\epsilon$ for $\lambda \succeq \lambda_0$, we find that for these $\lambda$,
\begin{equation}
\|X-X(1\otimes e_\lambda )\|_{A\tp A}<\epsilon .
\end{equation}
In particular, this is true for $X=(a\otimes 1)x$, $a\in A$, so that
$$
(a\otimes 1)\left[ x(1\otimes e_\lambda)-x\right] \longrightarrow 0.
$$

Next,
\begin{multline*}
x(1\otimes e_\lambda )(1\otimes a)=x(1\otimes (e_\lambda a))=x(1\otimes (a e_\lambda ))+ x(1\otimes (e_\lambda a-a e_\lambda))\\
=x(1\otimes a)(1\otimes e_\lambda )+x(1\otimes (e_\lambda a-a e_\lambda)).
\end{multline*}
The first summand tends to $x(1\otimes a)$ by (2.6), while the second summand tends to 0. The proofs of the remaining limit relations are similar. $\qquad \blacksquare$

\section{Multiplier Banach-Hopf algebras}

{\bf 3.1.} {\it Main notions}. Let $A$ be a polar Banach algebra over $\K$ with a nondegenerate product possessing a bounded approximate identity.

A nondegenerate continuous homomorphism $\Delta:\ A\to M(A\tp A)$ is called a {\it comultiplication} on $A$, if:
\begin{description}
\item[(i)] For any $a,b\in A$, the elements from $M(A\tp A)$ in the right-hand sides of the formulas
$$
T_1(a\otimes b)=\Delta (a)(1\otimes b);\quad T_2(a\otimes b)=(a\otimes 1)\Delta (b);
$$
$$
T_3(a\otimes b)=(1\otimes b)\Delta (a);\quad T_4(a\otimes b)=\Delta (b)(a\otimes 1)
$$
actually belong to $A\tp A$, so that $\Delta (a)\in M_0(A\tp A)$ for any $a\in A$.

\item[(ii)] The homomorphism $\Delta$ satisfies the following coassociativity condition -- the diagram
$$
\begin{CD}
A @>\Delta>> M(A\tp A) \\
@VV{\Delta}V @VV{\id \otimes \Delta}V \\
M(A\tp A) @>\Delta \otimes \id >> M(A\tp A\tp A)
\end{CD}
$$
in which the nondegenerate homomorphisms $\Delta$, $\Delta \otimes \id$, and $\id\otimes\Delta$, are extended to the appropriate multiplier algebras, -- is commutative.
\end{description}

A pair $(A,\Delta )$ is called a {\it regular multiplier Banach-Hopf algebra}, if the expressions $T_1,T_2,T_3$, and $T_4$ defined in (i) are extended to bijective linear isometric mappings $T_1,T_2,T_3,T_4:\ A\tp A \to A\tp A$. In this paper we will not consider non-regular algebras, thus the word ``regular'' will be dropped.

Let $\Delta'$ be the opposite comultiplication obtained from $\Delta$ by composing it with the flip $\sigma$ on $A\tp A$, $\sigma (a\otimes b)=b\otimes a$. Then $(A,\Delta')$ is also a multiplier Banach-Hopf algebra. Note that $\sigma$ is an isometry (\cite{BGR}, Section 2.1, Proposition 6(ii)) and can be extended to multipliers. For example, if $l\in L(A\tp A)$, then $l_\sigma =\sigma (l)$ acts as follows:
$$
l_\sigma (a\otimes b)=\sigma (l(\sigma (a\otimes b)))=\sigma (l(b\otimes a)),
$$
and the mapping $l\mapsto l_\sigma$ extends the flip operation on $A\otimes A$.

Let $m:\ A\tp A\to A$ be an extension by linearity and continuity of the multiplication operator, $m(a\otimes b)=ab$; see Proposition 2.1.7.1 in \cite{BGR}. Our construction of the counit follows \cite{VD94}. We define a continuous mapping $E:\ A\to L(A)$ setting
$$
E(a)b=mT_1^{-1}(a\otimes b).
$$
It is checked directly that $E(a)$ is indeed a left multiplier, and the identity
\begin{equation}
(\id \otimes E)((a\otimes 1)\Delta (b))=(ab)\otimes 1
\end{equation}
is valid for any $a,b\in A$.

The proof of (3.1) is identical with that of Lemma 3.2 in \cite{VD94}, with a single difference -- instead of a finite sum in the representation of $a\otimes b$ as an element from the range of $T_1$ we have to write a limit of an infinite sequence:
\begin{equation}
a\otimes b=\lim\limits_{j\to \infty} \sum\limits_{i=1}^{m_j}\Delta \left( a_i^{(j)}\right)\left( 1\otimes b_i^{(j)}\right)
\end{equation}
using the continuity of $\Delta$ and the polarity of $A$. In order to simplify notation, instead of the expressions like (3.2), we will always write
\begin{equation}
a\otimes b=\lim\sum\Delta (a_i)(1\otimes b_i),
\end{equation}
always remembering that (3.3) is just a shorthand for (3.2).

It follows from the surjectivity of $T_2$ that, for arbitrary elements $c,d\in A$,
$$
c\otimes d=T_2\left( \lim\sum a_j\otimes b_j\right) =\lim\sum (a_j\otimes 1)\Delta (b_j),
$$
and by the identity (3.1),
$$
(\id \otimes E)(c\otimes d)=\lim\sum (a_ib_i)\otimes 1
$$
where the existence of a limit $x=\lim\sum a_ib_i \in A\tp A$ follows from the continuity of $m$.

Thus $c\otimes E(d)=(\id \otimes E)(c\otimes d)=x\otimes 1$ for any $c,d\in A$, so that $E(A)\subset \K \cdot 1$.  Now we can define the {\it counit} $\varepsilon :\ A\to \K$ setting $\varepsilon (a)1=E(a)$.

\medskip
\begin{teo}
The counit $\varepsilon$ is a continuous homomorphism $A\to \K$, such that
\begin{gather}
\begin{split}
(\id\otimes \varepsilon)((a\otimes 1)\Delta (b)) & =ab;\quad (\varepsilon \otimes \id)(\Delta (a)(1\otimes b)) =ab;\\
(\id\otimes \varepsilon)(\Delta (a)(b\otimes 1)) & =ab;\quad (\varepsilon \otimes \id)((1\otimes a)\Delta (b)) =ab.
\end{split}
\end{gather}
\end{teo}

\medskip
For the {\it proof} see \cite{VD94}.$\qquad \blacksquare$

By the above construction, $\|\varepsilon \|\le 1$.

\medskip
The {\it antipode} $S:\ A\to M(A)$ is defined by the formula
\begin{equation}
S(a)b=(\varepsilon \otimes \id)T_1^{-1}(a\otimes b),\quad a,b\in A.
\end{equation}

\medskip
\begin{teo}
The antipode is a continuous anti-homomorphism $A\to M(A)$ satisfying the identities
\begin{equation}
m(S \otimes \id)(\Delta (a)(1\otimes b)) =\varepsilon (a)b;
\end{equation}
\begin{equation}
m(\id\otimes S)((a\otimes 1)\Delta (b)) =\varepsilon(b)a.
\end{equation}
The identities (3.6) and (3.7) define the antipode in a unique way. In addition, elements $S(a),a\in A$, belong actually to $A$.
\end{teo}

\medskip
{\it Proof}. Since $T_1$ is an isometry and $\|a\otimes b\|=\|a\| \cdot \|b\|$, we have
$$
\|S(a)b\| \le \|\varepsilon \otimes \id\|_{A\tp A\to A}\|a\|\cdot \|b\|,
$$
so that $\|S(a)\|\le C\|a\|$, which proves the continuity of $S$. The algebraic properties are proved just as in \cite{VD94} and \cite{Sun}. $\qquad\blacksquare$

Let $\varepsilon'$ and $S'$ be the counit and antipode corresponding to the opposite comultiplication $\Delta'$. Again, just as in the algebraic theory \cite{VD94}, we prove that $\varepsilon'=\varepsilon$, $S':\ A\to A$,
\begin{equation}
SS'=S'S=\id,
\end{equation}
and
\begin{equation}
(1\otimes Sb)\Delta (Sa)=(S\otimes S)(\Delta' (a)(1\otimes b)).
\end{equation}
The relation (3.8) also means that $S$ has a bounded inverse.

\medskip
\begin{prop}
The antipode $S$ can be extended to a continuous anti-homomorphism \linebreak $M(A)\to M(A)$.
\end{prop}

\medskip
{\it Proof}. We have to check that $S$ is nondegenerate, that is the set
$$
\operatorname{span}\{ S(a)b:\ a,b\in A\}
$$
is dense in $A$. Since $T_1$ is bijective, the set $\operatorname{span}\{ T_1^{-1}(a\otimes b):\ a,b\in A\}$ is dense in $A\tp A$. Let us choose $y\in A$ in such a way that $\varepsilon (y)=1$. Then $(\varepsilon \otimes \id)(y\otimes x)=x$ for every $x\in A$. This means that every element $x\in A$ can be approximated by linear combinations of elements $(\varepsilon \otimes \id)T_1^{-1}(a\otimes b)$ where the corresponding linear combinations of elements $a\otimes b$ are chosen to approximate $T_1(y\otimes x)$. $\qquad \blacksquare$

\bigskip
{\bf 3.2.} {\it Invariant functionals}. Let $\omega$ be a linear continuous mapping $A\to \K$, $a$ be an element from $A$. Define a left multiplier $\rho_l$ setting
$$
\rho_l(b)=(\omega \otimes \id)(\Delta (a)(1\otimes b))
$$
and a right multiplier
$$
\rho_r(b)=(\omega \otimes \id)((1\otimes b)\Delta (a)).
$$

The isometricity of $T_1$ and $T_3$ implies the continuity of $\rho_l$ and $\rho_r$, and it is checked easily that $(\rho_l,\rho_r)\in M(A)$. As in \cite{VD98}, we denote
$$
(\rho_l,\rho_r)=(\omega \otimes \id)\Delta (a).
$$
Similarly, using $T_2$ and $T_4$ we define $(\id \otimes \omega)\Delta (a)\in M(A)$.

The nonzero functional $\omega$ is called {\it left-invariant}, if
\begin{equation}
(\id \otimes \omega)\Delta (a)=\omega (a)1_{M(A)},
\end{equation}
and {\it right-invariant}, if
\begin{equation}
(\omega\otimes\id )\Delta (a)=\omega (a)1_{M(A)}.
\end{equation}

In this section and in Section 4, {\it we assume the existence of left-invariant and right-invariant functionals}. In fact, it suffices to know the existence of one of them -- if $\varphi$ is a left-invariant functional, and $\psi =\varphi \circ S$, then $\psi$ is right-invariant (see \cite{VD98}). In the examples of Section 5, the invariant functionals will be constructed explicitly.

Below we will prove the uniqueness of invariant functionals (up to a scalar factor). The above notation, $\varphi$ for a left-invariant functional, and $\psi$ for a right-invariant one, will be retained throughout the paper.

Let us study some properties of invariant functionals. First of all, they are faithful, that is the equality $\varphi (ba)=0$ (or the equality $\varphi (ab)=0$) for all $b\in A$ implies the equality $a=0$. The functional $\psi$ possesses a similar property. For the proof see \cite{VD98}, Proposition 3.4.

In the following lemma, often used in the sequel, the infinite sums (3.3) are understood in the sense of (3.2).

\medskip
\begin{lem}
Let $p_i,q_i\in A$. Suppose that there exists the limit $\lim\sum \Delta (p_j)(q_j\otimes 1)$ in $A\tp A$. Then, for any linear continuous functional $\chi: A\to \K$, there exists the limit $\lim\sum p_j\chi (q_j)$ in $A$.
\end{lem}

\medskip
{\it Proof}. We have $T_4(q_i\otimes p_i)=\Delta (p_i)(q_i\otimes 1)$,
$$
\sum p_i\otimes q_i=\sigma \circ T_4^{-1}\left( \sum \Delta (p_i)(q_i\otimes 1)\right).
$$
Therefore there exists $\lim\sum p_i\otimes q_i$, thus also $\lim\sum p_j\chi (q_j)$. $\qquad \blacksquare$

\medskip
Similar reasoning works for sums related to $T_1,T_2,T_3$.

Returning to properties of invariant functionals, we begin with the following result similar to Lemma 3.6 in \cite{VD98}. We consider only nonzero invariant functionals.

\medskip
\begin{lem}
The following sets of functionals on $A$ coincide:
\begin{equation}
\{ \psi (\cdot a):\ a\in A\}=\{ \varphi (\cdot a):\ a\in A\};
\end{equation}
\begin{equation}
\{ \psi (a \cdot ):\ a\in A\}=\{ \varphi (a\cdot ):\ a\in A\}.
\end{equation}
In particular, these sets do not depend on the choice of $\varphi,\psi$.
\end{lem}

\medskip
{\it Proof}. Let us prove the inclusion $\{ \psi (\cdot a):\ a\in A\}\subseteq\{ \varphi (\cdot a):\ a\in A\}$. Then the remaining three inclusions can be proved in a similar way.

In fact, we have to prove that for a given left-invariant functional $\varphi$, a given $a$, and for any right-invariant functional $\psi$, there exists such $c\in A$ that $\varphi (\cdot a)=\psi (\cdot c)$.

By definition, we have $(\id \otimes \varphi)\Delta (a)=\varphi (a)1_{M(A)}$. Choose $b\in A$ so that $\psi (b)=1$. Then, for any $x\in A$,
$$
\varphi (xa)=\psi (b)\varphi (xa)=(\psi \otimes \varphi )(\Delta (x)\Delta (a)(b\otimes 1)).
$$
Since $T_1$ is surjective, we may write
$$
\Delta (a)(b\otimes 1)=\lim\sum \Delta (c_i)(1\otimes d_i),\quad c_i,d_i\in A.
$$
Now
$$
\varphi (xa)=\lim\sum (\psi \otimes \varphi )(\Delta (x)\Delta (c_i)(1\otimes d_i))=\lim\sum \psi (xc_i)\varphi (d_i)=\psi (xc)
$$
where $c=\lim\sum c_i\varphi (d_i)$. The existence of this limit follows from Lemma 3.4. $\qquad \blacksquare$

\medskip
\begin{teo}
The left-invariant and right-invariant functionals are unique, up to multiplication by scalars.
\end{teo}

\medskip
The {\it proof}, based on Lemma 3.4 and Lemma 3.5, repeats the reasoning from \cite{VD98}.

\medskip
In fact, just as in \cite{VD98}, all the four sets of functionals listed in (3.12), (3.13) coincide (below this set of functionals will have the meaning of the dual object $\A$). The {\it proof}, identical with that of Proposition 3.11 in \cite{VD98}, uses also the following lemma needed later.

\medskip
\begin{lem}
Let $a,b,a_i,b_i\in A$. The following properties are equivalent:
\begin{description}
\item[(i)] $\Delta (a)(1\otimes b)=\lim\sum \Delta (a_i)(b_i\otimes 1)$;

\item[(ii)] $a\otimes S^{-1}b=\lim\sum (a_i\otimes 1)\Delta (b_i)$;

\item[(ii)] $(1\otimes a)\Delta (b)=\lim\sum (Sb_i)\otimes a_i$.
\end{description}
\end{lem}

\medskip
For the {\it proof} see \cite{VD94} (Lemma 5.5).

\bigskip
{\bf 3.3.} {\it Modular element. Modular automorphism}. A modular element $\delta$, a general version of the classical modular function defined on locally compact groups, is described in the next proposition.

\medskip
\begin{prop}
There exists such a multiplier $\delta \in M(A)$ that $(\varphi \otimes \id )\Delta (a)=\varphi (a)\delta$, that is
\begin{equation}
(\varphi \otimes \id )((1\otimes b)\Delta (a))=\varphi (a)b\delta;
\end{equation}
\begin{equation}
(\varphi \otimes \id )(\Delta (a)(1\otimes b))=\varphi (a)\delta b.
\end{equation}
\end{prop}

\medskip
{\it Proof}. For each $a\in A$, define a multiplier $\delta_a\in M(A)$ setting $\delta_a=(\varphi \otimes \id )\Delta (a)$, that is
$$
\delta_ab=(\varphi \otimes \id )(\Delta (a)(1\otimes b)),
$$
$$
b\delta_a=(\varphi \otimes \id )((1\otimes b)\Delta (a)).
$$
For a continuous linear functional $\omega$ on $A$, set
$$
\varphi_1(a)=(\varphi \otimes \omega )((1\otimes b)\Delta (a))=\omega (b\delta_a).
$$

Let us check that the functional $\varphi_1$ is left invariant. Since $A$ is polar, it is sufficient to prove that for any linear continuous functional $f$ on $A$ and any $c,d\in A$,
\begin{equation}
(f\otimes \varphi_1)(\Delta (c)(d\otimes 1))=\varphi_1(c)f(d);
\end{equation}
\begin{equation}
(f\otimes \varphi_1)((d\otimes 1)\Delta (c))=\varphi_1(c)f(d);
\end{equation}
see \cite{Sun}, Remark 6.1.3.

By Proposition 2.4, we can write $\Delta (c)$ as a strict limit
\begin{equation}
\Delta (c)=\lim\sum c_i'\otimes c_i'',\quad c_i',c_i''\in A.
\end{equation}
Writing the coassociativity relation for the opposite comultiplication $\Delta'$:
$$
(a\otimes 1\otimes 1)(\Delta'\otimes \id)(\Delta'(b)(1\otimes c))=(\id \otimes \Delta')((a\otimes 1)\Delta'(b))(1\otimes 1\otimes c)
$$
and applying the flip to both sides, we get the identity
\begin{equation}
(\Delta \otimes \id )((1\otimes a)\Delta (b))(c\otimes 1\otimes 1)=(1\otimes 1\otimes a)(\id \otimes \Delta)(\Delta (b)(c\otimes 1)),\quad a,b,c\in A.
\end{equation}

Using (3.18) and (3.19) we find that
\begin{multline*}
(f\otimes \varphi_1)(\Delta (c)(d\otimes 1))=\lim\sum f(c_i'd)\varphi_1(c_i'')=\lim\sum f(c_i'd)\omega (b\delta_{c_i''})\\
=\omega (\lim\sum (f\otimes \varphi \otimes 1)((c_i'd)\otimes ((1\otimes b)\Delta (c_i''))))=\omega ((f\otimes \varphi \otimes 1)((1\otimes 1\otimes b)(\id \otimes \Delta)(\Delta (c)(d\otimes 1))))\\
=\omega ((f\otimes \varphi \otimes 1)((\Delta \otimes \id)((1\otimes b)\Delta (c)))(d\otimes 1\otimes 1)).
\end{multline*}
Now we use (3.18) again, and then the left invariance of $\varphi$. Thus,
\begin{multline*}
(f\otimes \varphi_1)(\Delta (c)(d\otimes 1))=\lim\sum (f\otimes \varphi)(\Delta (c_i')(d\otimes 1))\omega (bc_i'')=\lim\sum \varphi (c_i')f(d)\omega (bc_i'')\\
=f(d)(\varphi \otimes \omega )((1\otimes b)\Delta (c))=f(d)\varphi_1(c),
\end{multline*}
so that we have proved (3.16).

In order to prove (3.17), we use the following coassociativity-like identity:
\begin{equation}
(1\otimes 1\otimes b)(\id \otimes \Delta )((c\otimes 1)\Delta (a))=(c\otimes 1\otimes 1)(\Delta \otimes \id)((1\otimes b)\Delta (a))
\end{equation}
obtained from the main identity $(\Delta \otimes \id)\Delta=(\id \otimes \Delta)\Delta$, where $\Delta:\ M(A)\to M(A\tp A)$ while $\Delta \otimes \id$ and $\id \otimes \Delta$ are extended to homomorphisms $M(A\tp A)\to M(A\tp A\tp A)$. Note that $(\Delta \otimes \id)(1\otimes b)=1\otimes 1\otimes b$, $(\id \otimes \Delta )(c\otimes 1)=c\otimes 1\otimes 1$.

Now, using (3.18) and (3.20) we get
\begin{multline*}
(f\otimes \varphi_1)((d\otimes 1)\Delta (c))=\lim\sum f(dc_i')\varphi_1(c_i'')=\lim\sum f(dc_i')\omega (b\delta_{c_i''})\\
=\omega \left( \lim\sum (f\otimes \varphi \otimes \id)((dc_i')\otimes (1\otimes b)\Delta (c_i''))\right)=\omega ((f\otimes \varphi \otimes \id)(1\otimes 1\otimes b)(\id \otimes \Delta)((d\otimes 1)\Delta (c)))\\
=\omega ((f\otimes \varphi \otimes \id)((d\otimes 1\otimes 1)(\Delta \otimes \id )((1\otimes b)\Delta (c))) =\lim\sum \omega ((f\otimes \varphi \otimes \id)((d\otimes 1)\Delta (c_i'))\otimes (bc_i''))\\
=\lim\sum (f\otimes \varphi )((d\otimes 1)\Delta (c_i'))\omega (bc_i'')=\lim\sum \varphi (c_i')f(d)\omega (bc_i'')=f(d)(\varphi \otimes \omega)((1\otimes b)\Delta (c))\\
=f(d)\varphi_1(c),
\end{multline*}
and we have proved (3.17).

By Theorem 3.6, $\varphi_1(a)=\lambda \varphi (a)$, where $\lambda \in \K$ does not depend on $a$. This means that
$$
\omega (b\delta_a)\varphi (c)=\omega (b\delta_c)\varphi (a)
$$
for any $a,b,c\in A$ and any functional $\omega$. Since $A$ is polar, this implies the identity $\varphi (c)\delta_a=\varphi (a)\delta_c$. Choose such an element $c$ that $\varphi (c)=1$ and set $\delta =\delta_c$. Then $\delta_a=\varphi (a)\delta$, so that $\delta$ is the required multiplier. $\qquad \blacksquare$

\medskip
The next proposition contains some identities for the modular element $\delta$.

\medskip
\begin{prop}
The modular element is invertible and satisfies the following equalities:
\begin{equation}
\Delta (\delta )=\delta \otimes \delta ;\quad \varepsilon (\delta )=1;\quad S(\delta )=\delta^{-1};
\end{equation}
\begin{equation}
\varphi (S(a))=\varphi (a\delta )\quad \text{for any $a\in A$}.
\end{equation}
In these formulas, $\Delta ,\varepsilon $, and $S$ are extended onto multipliers.
\end{prop}

\medskip
{\it Proof}. Applying $\Delta$ to both sides of the equality $\varphi (a)\delta =(\varphi \otimes \id )\Delta (a)$ we come to the identity
\begin{equation}
\varphi (a)\Delta (\delta )=(\varphi \otimes \id \otimes \id )((\id \otimes \Delta)\Delta (a)).
\end{equation}

Indeed, $(\varphi \otimes \id )\Delta (a)$ is the left multiplier
\begin{equation}
Lb=(\varphi \otimes \id )(\Delta (a)(1\otimes b)),\quad b\in A.
\end{equation}
Let us calculate the value of $\Delta (L)$ on the set of elements $\Delta (b)(v\otimes w)$, $v,w\in A$, whose span is dense in $A\tp A$. We have
$$
\Delta (L)(\Delta (b)(v\otimes w))=\Delta (Lb)(v\otimes w).
$$
Writing $\Delta (a)(1\otimes b)=\lim\sum a_i'\otimes a_i''$ we get
$$
Lb=\lim\sum \varphi(a_i')a_i'',\quad \Delta (Lb)=\lim\sum \varphi (a_i')\Delta (a_i''),
$$
so that
\begin{multline*}
\Delta (L)(\Delta (b)(v\otimes w))=\lim\sum \varphi (a_i')\Delta (a_i'')(v\otimes w)\\
=\lim\sum (\varphi \otimes \id \otimes \id)\left( a_i'\otimes (\Delta (a_i'')(v\otimes w))\right)=\lim\sum (\varphi \otimes \id \otimes \id )(((\id \otimes \Delta )(a_i'\otimes a_i''))(v\otimes w))\\
=(\varphi \otimes \id \otimes \id )((\id \otimes \Delta )(\Delta (a)(1\otimes b))(v\otimes w))=(\varphi \otimes \id \otimes \id )((\id \otimes \Delta )(\Delta (a))(\Delta (b)(v\otimes w)).
\end{multline*}

A similar equality holds for the appropriate right multipliers. Thus, (3.23) has been proved.

Next, in (3.23) we can use the coassociativity and obtain the equality
\begin{equation}
\varphi (a)\Delta (\delta )=(\varphi \otimes \id \otimes \id )((\Delta\otimes \id )\Delta (a)).
\end{equation}
Multiplying both sides from the right by the multiplier $1\otimes b$ and using again the definition of $\delta$ we find that
\begin{multline*}
(\varphi \otimes \id \otimes \id )((\Delta\otimes \id )\Delta (a))(1\otimes b)=\lim\sum (\varphi \otimes \id \otimes \id )\left( \Delta (a_i')\otimes a_i''\right)\\
=\lim\sum \varphi (a_i')\delta \otimes a_i''=\delta \otimes [(\varphi \otimes \id)(\Delta (a)(1\otimes b))] =\delta \otimes [\varphi (a)\delta b]
=\varphi (a)(\delta \otimes \delta )(1\otimes b),
\end{multline*}
and by (3.25),
$$
(\Delta (\delta )-\delta \otimes \delta )(1\otimes b)=0
$$
for any $b\in A$. Multiplying from the right by $c\otimes 1$, $c\in A$, we prove that $\Delta (\delta )-\delta \otimes \delta$ vanishes on a dense subset of $A\tp A$. This proves the first equality from (3.21).

Similarly, we calculate $\varepsilon (L)$ (for $L$ given in (3.24)) using (3.4):
$$
\varepsilon (L)\varepsilon (b)=\varepsilon (Lb)=\lim\sum \varphi (a_i')\varepsilon (a_i'')=(\varphi \otimes \varepsilon)(\Delta (a)(1\otimes b))=\varphi ((\id \otimes \varepsilon)(\Delta (a)(1\otimes b)))=\varphi (ab),
$$
so that
\begin{equation}
\varepsilon (L)\varepsilon (b)=\varphi (ab).
\end{equation}

Let us substitute for $b$ the approximate identity $e_\lambda$. We have $\varphi (ae_\lambda )\to \varphi (a)$, $\varepsilon (ze_\lambda)=\varepsilon (z)\varepsilon (e_\lambda)$ and $\varepsilon (ze_\lambda)\to \varepsilon (z)$ for each $z\in A$, so that $\varepsilon (e_\lambda)\to 1$. Passing to the limit in (3.26) we find that $\varepsilon (L)=\varphi (a)$.

On the other hand, it follows from the definition of $\delta$ that $\varepsilon (L)=\varphi (a)\varepsilon (\delta )$. Since $a$ is arbitrary, we get $\varepsilon (\delta )=1$.

To prove the identity for $S(\delta )$, we take $a,b\in A$ and write the identity
$$
\Delta (a\delta )(1\otimes b)=\Delta (a)(\delta \otimes (\delta b)).
$$
Let us apply $S\otimes \id$ to both sides and use the anti-homomorphism property of $S$. We obtain that
$$
(S\otimes \id )\Delta (a)(\delta \otimes (\delta b))=(S(\delta )\otimes 1)(S\otimes \id )(\Delta (a)(1\otimes \delta b).
$$
Taking an arbitrary $c\in A$ and using (3.7) we get the identity
\begin{equation}
c\varepsilon (a\delta )b=S(\delta)c\varepsilon (a)\delta b.
\end{equation}

In (3.27), we specify $c=e_\lambda$ and pass to the limit using the fact that $\varepsilon (\delta )=1$. This results in the equality $b=S(\delta )\delta b$, so that $S(\delta )\delta =1$.

Similarly we can write
$$
(c\otimes 1)\Delta (\delta a)=(1\otimes \delta )(c\delta \otimes 1)\Delta (a)
$$
whence
$$
(\id \otimes S)((c\otimes 1)\Delta (\delta a)(1\otimes b))=(\id \otimes S)(((c\delta)\otimes 1)\Delta (a)(1\otimes b))(1\otimes S(\delta )).
$$
Applying the mapping $m$ and using the identity (3.6) we find that $c\varepsilon (\delta a)b=c\delta \varepsilon (a)bS(\delta )$, so that $cb=c\delta bS(\delta )$. Here we set $b=e_\lambda$, pass to the limit and use the arbitrariness of $c$. As a result, $\delta S(\delta )=1$. This means that $\delta$ is invertible and $S(\delta )=\delta^{-1}$.

The proof of the identity (3.22) is similar to that of Proposition 3.10 in \cite{VD98}. $\qquad\blacksquare$

\medskip
Let $\mathfrak S$ be a topology on $A$ generated by the set of seminorms $x\mapsto |\Phi (x)|$ where $\Phi$ is a functional from the family (3.12)-(3.13).

\medskip
\begin{prop}
There exists an $\mathfrak S$-continuous automorphism $\beta$ of the algebra $A$, such that
\begin{equation}
\varphi (ab)=\varphi (b\beta (a))\quad \text{for all $a,b\in A$}.
\end{equation}
The functional $\varphi$ is $\beta$-invariant.
\end{prop}

\medskip
{\it Proof}. Note first of all that $\operatorname{span}\{ab:\ a,b\in A \}$ is dense in $A$. Indeed, we assume that the mapping $T_1(a\otimes b)=\Delta (a)(1\otimes b)$ extends to a bijection of $A\tp A$ onto itself. On the other hand,
$$
(\varepsilon \otimes \id )(\Delta (a)(1\otimes b))=ab.
$$

Choose $c\in A$ in such a way that $\varepsilon (c)=1$. For any $d\in A$, $(\varepsilon \otimes \id )(c\otimes d)=d$. Since $T_1$ is a bijection, we can write $c\otimes d=\lim\sum \Delta (a_i)(1\otimes b_i)$, $d=\lim\sum a_ib_i$.

Define $\beta$ by the relation (3.28); that is possible since $\varphi$ is faithful. For any $a,b,c\in A$,
$$
\varphi (abc)=\varphi (c\beta (ab))=\varphi (bc\beta (a))=\varphi (c\beta (a)\beta (b)),
$$
and the faithfulness of $\varphi$ implies the equality $\beta (ab)=\beta (a)\beta (b)$. Similarly, the linearity of $\varphi$ implies the linearity of $\beta$.

If $\beta (a_1)=\beta (a_2)$, $a_1,a_2\in A$, then
$$
\varphi ((a_1-a_2)b)=\varphi (b(\beta (a_1)-\beta (a_2))=0
$$
for any $b\in A$, thus $a_1=a_2$. The surjectivity of $\beta$ follows from the above coincidence of the families of functionals (3.12) and (3.13).

Let us prove the $\mathfrak S$-continuity of $\beta$. Let $a_\lambda\to a$ be a $\mathfrak S$-convergent net in $A$. Then $\varphi (a_\lambda b)\to \varphi (ab)$ for any $b\in A$, so that
$$
\varphi (b[\beta (a_\lambda )-\beta (a)])=\varphi (a_\lambda b)-\varphi (ab)\to 0,
$$
which means the $\mathfrak S$-continuity of $\beta$.

To prove that $\varphi$ is $\beta$-invariant, we write (3.28) with $b=e_\lambda$ and pass to the limit. $\qquad\blacksquare$

\medskip
As in \cite{VD98} (page 340), we get also the relation
\begin{equation}
\varphi (\delta^{-1}a\delta )=\tau \varphi (a)\quad a\in A,
\end{equation}
with some constant $\tau \in \K$.

Let $\beta'$ be a similar automorphism associated with the right-invariant functional $\psi$. Then $S\beta'=\beta^{-1}S$. For the proof see \cite{VD98}, Proposition 3.13.

Following \cite{VD98} (Propositions 3.14 and 3.15) and \cite{Ti}, Section 2.2.4, we prove, for any $a\in A$, the relations
$$
\Delta (\beta (a))=(S^2\otimes \beta )\Delta (a); \quad \Delta (\beta' (a))=(\beta'\otimes S^{-2})\Delta (a);
$$
$$
\beta (\delta )=\beta' (\delta )=\frac1\tau \delta;\quad \delta\beta (a)=\beta'(a)\delta
$$
where $\tau \in \K$ is the constant appearing in (3.29).

\section{The dual object}

{\bf 4.1.} {\it The dual algebra}. Denote by $\A$ the set of linear continuous functionals on $A$ of the form $\varphi (\cdot a)$, $a\in A$. Below we assume that the left-invariant functional $\varphi$ is such that $\| \varphi \|\le 1$ and possesses the {\it norm reproducing property}
\begin{equation}
\|a\|=\sup\limits_{x\ne 0}\frac{|\varphi (xa)|}{\|x\|},\quad a\in A.
\end{equation}

Equivalently, functionals from $\A$ can be represented as $\varphi (b\cdot)$, $\psi (\cdot c)$, $\psi (d\cdot )$ where $b,c,d\in A$, $\psi$ is a right-invariant functional.

$\A$ is obviously a Banach space over $\K$ with respect to the norm $\|\varphi (\cdot a)\|=\|a\|$, well-defined since $\varphi$ is faithful. By (4.1), this norm coincides with the standard norm of the functional $\varphi (\cdot a)$ on $A$. The functional $\varphi (\cdot a)$ can be seen as the Fourier transform of an element a. Then the identity (4.1) can be interpreted as a kind of the Plancherel formula.

The product in $\A$ is defined as follows. If $\omega_1=\varphi (\cdot a_1)$, $\omega_2=\varphi (\cdot a_2)$, $a_1,a_2\in A$, then, by definition,
\begin{equation}
(\omega_1\omega_2)(x)=(\omega_1\otimes \omega_2)\Delta (x)
\end{equation}
where
\begin{equation}
(\omega_1\otimes \omega_2)\Delta (x)=(\varphi \otimes \varphi )(\Delta (x)(a_1\otimes a_2)).
\end{equation}
The right-hand side of (4.3) is well-defined for all $x\in A$; it defines a linear continuous functional.

More specifically, writing $a_1\otimes a_2=\lim\sum \Delta (p_i)(q_i\otimes 1)$, $p_i,q_i\in A$, we have
\begin{equation}
(\omega_1\omega_2)(x)=\lim\sum \varphi (xp_i)\varphi (q_i)=\varphi (xb)
\end{equation}
where $b=\lim\sum p_i\varphi (q_i)\in A$ exists by virtue of Lemma 3.4. Therefore $\omega_1\omega_2\in \A$.

This product is nondegenerate (see \cite{VD98}, page 346) and associative (see the proof in \cite{VD94} for more general functionals); this associativity is based on the coassociativity of $\Delta$.

Looking at the element $b$ in (4.4) we find, as in the proof of Lemma 3.4, that
\begin{multline*}
b=(\id \otimes \varphi )(\lim\sum p_i\otimes q_i)=(\id \otimes \varphi )(\sigma \circ T_4^{-1}(\lim\sum \Delta (p_i)(q_i\otimes 1)))\\
=(\id \otimes \varphi )(\sigma \circ T_4^{-1}(a_1\otimes a_2)).
\end{multline*}
Since $\|\varphi \|\le 1$ and $T_4$ is an isometry, we find that $\|b\|\le \|a_1\otimes a_2\|=\|a_1\|\cdot \|a_2\|$, so that $\|\omega_1\omega_2\|\le \|\omega_1\|\cdot \|\omega_2\|$, and $\A$ is indeed a Banach algebra.

\bigskip
{\bf 4.2.} {\it Comultiplication on $\A$}. In order to define $\D:\ \A\to M(\A \tp \A)$, we first define, for any $\omega_1,\omega_2\in \A$, the linear continuous functionals $\D (\omega_1)(1\otimes \omega_2)$ and $(\omega_1\otimes 1)\D (\omega_2)$ on $A\tp A$ as follows:
\begin{equation}
(\D (\omega_1)(1\otimes \omega_2))(x\otimes y)=(\omega_1\otimes \omega_2)((x\otimes 1)\Delta (y)),
\end{equation}
\begin{equation}
((\omega_1\otimes 1)\D (\omega_2))(x\otimes y)=(\omega_1\otimes \omega_2)(\Delta (x)(1\otimes y)),
\end{equation}
$x,y\in A$.

\medskip
\begin{lem}
If $\omega_1,\omega_2,\omega_3\in \A$, then the functionals $\D (\omega_2)(1\otimes \omega_3)$ and $(\omega_1\otimes 1)\D (\omega_2)$ defined in (4.5), (4.6) belong to $\A \tp \A$, and
\begin{equation}
((\omega_1\otimes 1)\D (\omega_2))(1\otimes \omega_3)=(\omega_1\otimes 1)(\D (\omega_2)(1\otimes \omega_3)).
\end{equation}
\end{lem}

\medskip
{\it Proof}. Let $\omega_j=\varphi (\cdot a_j)$, $j=2,3$. Writing $a_2\otimes a_3=\lim\sum \Delta (b_i)(c_i\otimes 1)$, $b_i,c_i\in A$, we find that
$$
(\D (\omega_2)(1\otimes \omega_3))(x\otimes y)=(\varphi \otimes \varphi )((x\otimes 1)\Delta (y)(a_2\otimes a_3))=\lim\sum \varphi (xc_i)\varphi (yb_i).
$$
Thus,
\begin{equation}
\D (\omega_2)(1\otimes \omega_3)=\lim\sum \varphi (\cdot c_i)\varphi (\cdot b_i)
\end{equation}
(the convergence of functionals in the sense of $\A \tp \A$ is proved as in Lemma 3.4).

Similarly, using the representation $\omega_j=\psi (a_j'\cdot )$, $j=1,2$, we get
\begin{equation}
(\omega_1\otimes 1)\D (\omega_2)=\lim\sum \psi (c_j'\cdot )\psi (b_j'\cdot ).
\end{equation}

Let us write (4.8) as a convergent limit relation for functionals:
$$
\D (\omega_2)(1\otimes \omega_3)=\lim\sum\xi_i'\otimes \xi_i'',\quad \xi_i',\xi_i''\in \A.
$$
We have
\begin{multline*}
((\omega_1\otimes 1)(\D (\omega_2)(1\otimes \omega_3))) (x\otimes y)=\lim\sum ((\omega_1\xi_i')\otimes \xi_i'')(x\otimes y)\\
=\lim\sum ((\omega_1\otimes \xi_i')\Delta (x))\xi_i''(y)=\lim\sum ((\omega_1\otimes \xi_i'\otimes \xi_i'')(\Delta (x)\otimes y).
\end{multline*}

Let us approximate $\Delta (x)=\lim\sum \mu_j'\otimes \mu_j''$ in the topology of $M_0(A\tp A)$. The convergence in $i$ is uniform with respect to $j$ because the convergence in $M_0(A\tp A)$ means the convergence of multipliers in strong operator topology, which implies their uniform boundedness (\cite{PGS}, Theorem 2.1.20). Therefore we may change the order of convergence:
\begin{multline*}
((\omega_1\otimes 1)(\D (\omega_2)(1\otimes \omega_3))) (x\otimes y)=\lim\limits_j\sum\lim\limits_i\sum \omega_1(\mu_j')(\xi_i'\otimes \xi_i'')(\mu_j''\otimes y)\\
=\lim\sum \omega_1(\mu_j')(\D (\omega_2)(1\otimes \omega_3))(\mu_j''\otimes y)=\lim\sum \omega_1(\mu_j')(\omega_2\otimes \omega_3)((\mu_j''\otimes 1)\Delta (y)).
\end{multline*}

As above, we approximate $\Delta (y)=\lim\sum \nu_k'\otimes \nu_k''$ in $M_0(A\tp A)$ and get
\begin{equation}
((\omega_1\otimes 1)(\D (\omega_2)(1\otimes \omega_3)))(x\otimes y)=\lim\sum\lim\sum \omega_1(\mu_j')\omega_2(\mu_j''\nu_k')\omega_3(\nu_k'').
\end{equation}

In a similar way, we transform the left-hand side of (4.7) and obtain for it the expression identical with the one in the right-hand side of (4.10). This proves the identity (4.7). $\qquad \blacksquare$

\medskip
The identity (4.7) means that the mappings $\omega_1\otimes 1\mapsto (\omega_1\otimes 1)\D (\omega_2)=\Hat{T}_2(\omega_1\otimes \omega_2)$ and $1\otimes \omega_3\mapsto \D (\omega_2)(1\otimes \omega_3)=\Hat{T}_1(\omega_2\otimes \omega_3)$ extend to the right multiplier $\D (\omega_2)_r$ and the left multiplier $\D (\omega_2)_l$ respectively, and this pair defines a multiplier $\D (\omega_2)\in M(\A \tp \A)$.

By (4.4), $\Hat{T}_1(\omega_1\otimes \omega_2)(x\otimes y)=(\omega_1\otimes \omega_2)(T_2(x\otimes y))$, for any $x,y\in A$, $\omega_1,\omega_2\in \A$. As before, we write $\omega_1=\varphi (\cdot a_1)$, $\omega_2=\varphi (\cdot a_2)$, but here it will be convenient to represent
$$
\Delta (a_2)(a_1\otimes 1)=\lim\sum p_i\otimes q_i,\quad p_i,q_i\in A.
$$
The limit is in the topology of $A\tp A$. Then, for any $x,y\in A$,
\begin{multline*}
((\omega_1\otimes \omega_2)(x\otimes y)=\varphi (xa_1)\varphi (ya_2)=(\varphi \otimes \varphi )((x\otimes 1)\Delta (ya_2)(a_1\otimes 1))\\
=\lim\sum (\varphi \otimes \varphi )((x\otimes 1)\Delta (y)(p_i\otimes q_i)).
\end{multline*}
Writing $\omega_i'=\varphi (\cdot p_i)$, $\omega_i''=\varphi (\cdot q_i)$, we find that
$$
(\omega_1\otimes \omega_2)(x\otimes y)=\lim\sum (\omega_i'\otimes \omega_i'')((x\otimes 1)\Delta (y))=\lim\sum (\D (\omega_i')(1\otimes \omega_i''))((x\otimes y),
$$
where the limit is in the sense of $\A \tp \A$. This proves the surjectivity of $\Hat{T}_1$.

Extending by linearity and continuity we find that
$$
\Hat{T}_1(\omega_1\otimes \omega_2)(p)=(\omega_1\otimes \omega_2)(T_2p)
$$
for all $p\in \A \tp \A$. Since $T_2$ is an isometry, we have
$$
\left| \Hat{T}_1(\omega_1\otimes \omega_2)(p)\right| \le \|\omega_1\otimes \omega_2\|_{\A \tp \A}\cdot \|p\|,
$$
so that $\left\| \Hat{T}_1(\omega_1\otimes \omega_2)\right\|_{\A \tp \A}\le \|\omega_1\otimes \omega_2\|_{\A \tp \A}$.

Similarly we write $(\omega_1\otimes \omega_2)(p)=\Hat{T}_1(\omega_1\otimes \omega_2)(T_2^{-1}p)$ and obtain the inverse inequality. Therefore $\Hat{T}_1$ {\it is an isometry}.
Note the importance of the norm reproducing property (4.1). Here we interpreted the norm on $\A$ as the standard norm of functionals, while the proof that $\A$ is a Banach algebra was based on the fact that the norm of the functional $\varphi (\cdot a)$ equals $\|a\|$.

In a similar way we define the mappings $\Hat{T}_2,\Hat{T}_3$, and $\Hat{T}_4$, and prove their isometry and surjectivity properties. In fact, we followed \cite{Sun} (pages 90, 91). The difference from the purely algebraic case is the need to check the possibility to change the order of limits, and that is done as in the proof of Lemma 4.1.

\medskip
\begin{prop}
The mapping $\D$ is a continuous homomorphism $\A\to M(\A \tp \A)$ satisfying the coassociativity condition. The counit $\Hat \varepsilon:\ \A\to \K$ and the antipode $\Hat S:\A\to \A$ given by the formulas
\begin{equation}
\Hat \varepsilon (\varphi (\cdot a))=\Hat \varepsilon (\varphi (a\cdot ))=\varphi (a),
\end{equation}
\begin{equation}
\Hat \varepsilon (\psi (\cdot a))=\Hat \varepsilon (\psi (a\cdot ))=\psi (a),
\end{equation}
\begin{equation}
\Hat S(\omega )(a)=\omega (S(a)),\quad a\in A,\omega \in \A,
\end{equation}
are continuous homomorphisms satisfying, for any $\omega_1,\omega_2\in \A$, the identities
\begin{equation}
(\Hat \varepsilon \otimes \id)(\D (\omega_1)(1\otimes \omega_2))=\omega_1\omega_2
\end{equation}
and three other identities similar to (3.4),
\begin{equation}
m(\Hat S\otimes \id)(\D (\omega_1)(1\otimes \omega_2))=\Hat \varepsilon (\omega_1)\omega_2,
\end{equation}
\begin{equation}
m(\id\otimes \Hat S)((\omega_1\otimes 1)\D (\omega_2))=\Hat \varepsilon (\omega_2)\omega_1.
\end{equation}
\end{prop}

\medskip
{\it Proof}. The proof of the statement about $\D$ is similar to that of Proposition 4.6 from \cite{VD98} and Proposition 7.1.7 from \cite{Sun}.

Let $\omega \in \A$. Then we can write
$$
\omega =\varphi (a\cdot )=\varphi (\cdot a')=\psi (b\cdot )=\psi (\cdot b')
$$
for some $a,a',b,b'\in A$. In order to justify the definitions (4.11)-(4.12), we have to check that $\varphi (a)=\varphi (a')=\psi (b)=\psi (b')$. Since $\varphi (ax)=\varphi (x\beta (a))$ (see (3.28)) and $\varphi =\varphi \circ \beta$, we find that $a'=\beta (a)$ and $\varphi (a')=\varphi (a)$. In order to compare $\varphi (a)$ and $\psi (b)$, we write $\varphi (ax)=\psi (bx)$, set $x=e_\lambda$, and pass to the limit.

Let us prove (4.14). Writing $\omega_1=\varphi (\cdot a_1)$, $\omega_2=\varphi (\cdot a_2)$, we get
$$
(\D (\omega_1)(1\otimes \omega_2))(x\otimes y)=(\omega_1\otimes \omega_2)((x\otimes 1)\Delta (y))=(\varphi \otimes \varphi)((x\otimes 1)\Delta (y)(a_1\otimes a_2)).
$$
Let $a_1\otimes a_2=\lim\sum \Delta (p_i)(q_i\otimes 1)$. Then
\begin{equation}
(\D (\omega_1)(1\otimes \omega_2))(x\otimes y)=\lim\sum \varphi (xq_i)\varphi (yp_i)
\end{equation}
due to the left invariance of $\varphi$.

On the other hand,
\begin{equation}
(\omega_1\omega_2)(y)=(\varphi\otimes \varphi )(\Delta (y)(a_1\otimes a_2))=\lim\sum \varphi (q_i)\varphi (yp_i).
\end{equation}
Applying $\Hat \varepsilon \otimes \id$ to both sides of (4.17) and comparing with (4.18), we prove (4.14). The proofs of other identities for $\Hat \varepsilon$ are similar.

Note that (4.14) implies the homomorphism property of $\Hat \varepsilon$. Indeed, let us take $\omega_3\in \A$ and write (4.14) as
\begin{equation}
(\Hat \varepsilon \otimes \id) (\D (\omega_1\omega_3)(1\otimes \omega_2))=\omega_1\omega_3\omega_2.
\end{equation}
The left-hand side of (4.19) equals
$$
(\Hat \varepsilon \otimes \id) (\D (\omega_1)(\D (\omega_3)(1\otimes \omega_2)),
$$
while the right-hand side is equal to
$$
\omega_1(\omega_3\omega_2)=\omega_1(\Hat \varepsilon \otimes \id)(\D (\omega_3)(1\otimes \omega_2)).
$$
Due to the surjectivity of $\Hat{T}_1$, we may substitute in both sides $\omega_3\otimes \omega_2$ for $\D (\omega_3)(1\otimes \omega_2)$. We obtain that
$$
(\Hat \varepsilon \otimes \id) (\D (\omega_1)(\omega_3\otimes \omega_2))=\omega_1(\Hat \varepsilon \otimes \id)(\omega_3\otimes \omega_2)=\Hat \varepsilon (\omega_3)\omega_1\omega_2,
$$
so that
$$
(\Hat \varepsilon \otimes \id) (\D (\omega_1)(1\otimes \omega_2)(\omega_3\otimes 1))=\Hat \varepsilon (\omega_3)(\Hat \varepsilon \otimes \id) (\D (\omega_1)(1\otimes \omega_2)).
$$

Here we substitute $\omega_1\otimes \omega_2$ in both sides for $\D (\omega_1)(1\otimes \omega_2)$. This results in the equality
$$
(\Hat \varepsilon \otimes \id)((\omega_1\omega_3)\otimes \omega_2)=\Hat \varepsilon (\omega_3)(\Hat \varepsilon \otimes \id)(\omega_1\otimes \omega_2),
$$
thus $\Hat \varepsilon (\omega_1\omega_3)\omega_2=\Hat \varepsilon (\omega_1)\Hat \varepsilon (\omega_3)\omega_2$, and since the algebra $\A$ is nondegenerate, $\Hat \varepsilon (\omega_1\omega_3)$ \linebreak
$=\Hat \varepsilon (\omega_1)\Hat \varepsilon (\omega_3)$.

Turning to the antipode $\Hat S$ on $\A$ defined by (4.13), we use the standard identities for the antipode $S$ on $A$ extended onto $M(A)$, see \cite{VDW}, Section 2:
\begin{equation}
T_1^{-1}(a\otimes b)=((\id \otimes S)\Delta (a))(1\otimes b);
\end{equation}
\begin{equation}
T_2^{-1}(a\otimes b)=(a\otimes 1)((S\otimes\id )\Delta (b)),
\end{equation}
$a,b\in A$.

Let $\omega_1,\omega_2\in \A$, $\omega_1=\varphi (a_1\cdot )$, $a_1\in A$. Then by (4.5) and (4.21), for any $x\in A$,
\begin{multline*}
(m(\Hat S\otimes \id)(\D (\omega_1)(1\otimes \omega_2)))(x)=(\D (\omega_1)(1\otimes \omega_2))((S\otimes \id)\Delta (x))\\
=(\omega_1\otimes \omega_2)(T_2((S\otimes \id)\Delta (x)))=(\varphi \otimes \omega_2)((a_1\otimes 1)T_2((S\otimes \id)\Delta (x))\\
=(\varphi \otimes \omega_2)(T_2((a_1\otimes 1)(S\otimes \id)\Delta (x)))=(\varphi \otimes \omega_2)(a_1\otimes x)=\varphi (a_1)\omega_2 (x)\\
=\Hat \varepsilon (\omega_1)\omega_2 (x),
\end{multline*}
so that we have obtained the first defining identity (4.15) for the antipode.

Similarly, using (4.6) and (4.20), we prove the second defining identity (4.16). As we mentioned in Theorem 3.2, the identities (4.15) and (4.16) are sufficient to define the antipode.

The continuity properties of the mappings treated in this proposition are obvious. $\qquad \blacksquare$

\bigskip
{\bf 4.3.} {\it Invariant functionals}. The right-invariant continuous linear functional $\Hat \psi$ and the left-invariant continuous linear functional $\Hat \varphi$ on $\A$ are given by the formulas
$$
\Hat \psi (\omega )=\varepsilon (a)\quad \text{for $\omega =\varphi (\cdot a)$};
$$
$$
\Hat \varphi (\omega )=\varepsilon (a)\quad \text{for $\omega =\psi (a\cdot )$};
$$
The proof is similar to that in \cite{VD98} (Proposition 4.8).

The following important lemma is proved just like its algebraic counterpart (\cite{VD98}, Lemma 4.11).

\medskip
\begin{lem}
Let $\omega ,\omega_1\in \A$, $\omega =\varphi (\cdot a)$. Then $\Hat \psi (\omega_1\omega)=\omega_1(S^{-1}(a))$.
\end{lem}

\bigskip
{\bf 4.4.} {\it Biduality theorem}. Details of the proof of the following biduality property, fundamental for the duality theory, were kindly provided to the author by L. Vainerman.

For any $a\in A$, define a functional $\Gamma (a):\ \A\to \K$ setting
$$
\Gamma (a)(\omega )=\omega (a),\quad \omega\in \A.
$$
By Lemma 4.3, if we write $\omega =\varphi (\cdot S(a))$, then $\Gamma (a)=\Hat \psi (\cdot \omega )$, so that $\Gamma (a)\in \AAA$ where $\AAA$ denotes the dual object to $\A$. Obviously, the mapping $a\mapsto \Gamma (a)$ is a continuous isomorphism of Banach spaces.

Below we assume that the Banach algebra $\A$ possesses a bounded approximate identity, so that the results on extension of homomorphism onto multipliers are applicable. We also assume the norm reproducing property for $\A$. It would be interesting to find some sufficient conditions for these properties formulated in terms of $A$.

\medskip
\begin{teo}
$\Gamma$ is an isomorphism of multiplier Banach-Hopf algebras, that is
\begin{equation}
\Gamma (a_1\cdot a_2)=\Gamma (a_1)\cdot \Gamma (a_2),
\end{equation}
for any $a_1,a_2\in A$, and
\begin{equation}
\Hat{\D}(\Gamma (a))=(\Gamma \otimes \Gamma)(\Delta (a)),
\end{equation}
for any $a\in A$. Here $\Hat{\D}$ is the coproduct map in $\AAA$.
\end{teo}

\medskip
{\it Proof}. First, one can equivalently define $\Gamma$ as follows:
\begin{equation}
\Gamma (a)(\omega )=\Hat \psi (\omega \cdot \widehat{S(a)}),
\end{equation}
for all $a\in A$, $\omega \in \A$; here $\widehat{S(a)}=\varphi (\cdot S(a))$ is the Fourier transform of $S(a)$. Indeed, the right-hand side of (4.24) is equal to $\omega (S^{-1}(S(a)))=\omega (a)$ due to Lemma 4.3.

Then we find, using Lemma 4.3 again, that
\begin{multline*}
(\Gamma (a_1)\cdot \Gamma (a_2))(\omega )=(\Gamma (a_1)\otimes \Gamma (a_2))(\D (\omega ))=(\Hat \psi \otimes \Hat\psi )(\D (\omega)(\widehat{S(a_1)}\otimes \widehat{S(a_2)}))\\
=\Hat \psi \left( \left[ (\id \otimes \Hat \psi)(\D (\omega )(1\otimes \widehat{S(a_2)}))\right]\cdot \widehat{S(a_1)}\right)=\left[ (\id \otimes \Hat \psi)(\D (\omega )(1\otimes \widehat{S(a_2)}))\right](a_1).
\end{multline*}
Thus, in order to prove the needed equality for $\Gamma$, it suffices to show that the functional in the square brackets is equal to $\omega (\cdot a_2)$. Now, we can choose $\omega$ in the form of $\Hat b=\varphi (\cdot b)$. Then everything we need can be formulated as the following equality:
\begin{equation}
(\id \otimes \Hat \psi)(\D (\Hat b)(1\otimes \widehat{S(a)}))=\widehat{ab},
\end{equation}
for any $a,b\in A$.

To prove (4.25), we use the reasoning from the proof of Proposition 4.8 in \cite{VD98} -- write $b\otimes S(a)=\lim\sum \Delta (p_i)(q_i\otimes 1)$, then the calculation in the beginning of that proof means that
$$
\D (\Hat b)(1\otimes \widehat{S(a)})=\lim\sum \left( \widehat{q_i}\otimes \widehat{p_i}\right) ,
$$
so that
$$
(\id \otimes \Hat \psi)(\D (\Hat b)(1\otimes \widehat{S(a)}))=\lim\sum \Hat\psi (\widehat{p_i})\widehat{q_i}=\lim\sum \varepsilon (p_i)\widehat{q_i}.
$$

It remains to show that the last expression equals $ab$. Using properties of the antipode we can write the above equality for $b\otimes S(a)$ as
$$
a\otimes b=\lim\sum (S^{-1}\otimes \id )(\sigma \Delta (p_i))(1\otimes q_i).
$$
Looking at this equality as an equality in the opposite multiplier Banach-Hopf algebra whose coproduct is $\sigma \Delta$, the antipode is $S^{-1}$, the multiplication and the counit being the same, applying to both its sides the multiplication map $m$ and using the identity (3.6), we finish the proof of (4.22).

The equality (4.23) is an equality in $M(\AAA \tp \AAA)$, and the map in the right-hand side is understood as follows. We have seen that the map $\Gamma:\ A\to \AAA$ is multiplicative, so that the map $\Gamma \otimes \Gamma:\ A\tp A\to \AAA \tp \AAA$ is a well-defined multiplicative map. It can be seen as a multiplicative map from $A\tp A$ to $M(\AAA \tp \AAA)$, which has a canonical extension to a multiplicative map from $M(A\tp A)$ to $M(\AAA \tp \AAA)$. This extension is exactly the map in the right-hand side of (4.23).

As we know, the multiplier $\Hat{\D} (\Gamma (a))$ is defined by the elements $\Hat{\D} (\Gamma (a))(1\otimes \Gamma (b))$ and $(\Gamma (b)\otimes 1)\Hat{\D} (\Gamma (a))$ of $\AAA\tp \AAA$. The equality of multipliers in (4.23) is equivalent to the system of the following equalities on $\AAA\tp \AAA$:
\begin{equation}
\Hat{\D} (\Gamma (a))(1\otimes \Gamma (b))=(\Gamma \otimes \Gamma)(\Delta (a))(1\otimes \Gamma (b));
\end{equation}
\begin{equation}
(\Gamma (b)\otimes 1)\Hat{\D} (\Gamma (a))=(\Gamma (b)\otimes 1)(\Gamma \otimes \Gamma)(\Delta (a)),
\end{equation}
for any $a,b\in A$. Let us prove (4.26); the proof of (4.27) is similar.

The right-hand side of (4.26) evaluated on $\omega_1\otimes \omega_2\in \A\otimes \A$ gives $(\omega_1\otimes \omega_2)(\Delta (a)(1\otimes b))$ by the definition of $\Gamma$ and by its multiplicativity. The left-hand side of (4.26) evaluated on $\omega_1\otimes \omega_2$ can be rewritten, using (4.5), as $(\Gamma (a)\otimes \Gamma (b))(\omega_1\otimes 1)\D (\omega_2))$, which is equal, by the definition of $\Gamma$, to $((\omega_1\otimes 1)\D (\omega_2))(a\otimes b)$. Finally, the identity (4.6) shows that this expression is equal to $(\omega_1\otimes \omega_2)(\Delta (a)(1\otimes b))$. This proves (4.26). $\qquad \blacksquare$

\bigskip
{\bf 4.5.} {\it Multiplier Banach-Hopf algebras of compact and discrete types}. A multiplier Banach-Hopf algebra $(A,\Delta )$ with invariant functionals is said to be of compact type, if $A$ is a unital algebra, and to be of discrete type, if there is such a nonzero element $h\in A$ that $ah=\varepsilon (a)h$ for all $a\in A$.

\medskip
\begin{teo}
If $(A,\Delta )$ is of discrete type, then the dual object $(\A ,\D )$ is of compact type. If $(A,\Delta )$ is of compact type, then $(\A ,\D )$ is of discrete type.
\end{teo}

\medskip
The {\it proof} is identical to the one given in \cite{VD98} (Proposition 5.3) for the purely algebraic case.

\section{Examples}

In this section, we describe the three examples listed in Introduction. For each case, we give explicit descriptions of the multiplication, comultiplication, counit and antipode, invariant functionals, identify the bounded approximate identities and prove the norm reproducing property for the initial and dual algebras.

\bigskip
{\bf 5.1.} {\it Discrete groups}. Let $G_1$ be a discrete group. Let $A_1=c_0(G_1)$ be the commutative Banach algebra of $\K$-valued functions on $G_1$ tending to zero by the filter of complements to finite sets, with the pointwise operations and $\sup$-norm.

\medskip
\begin{lem}
The multiplier algebra $M(A_1)$ is isomorphic to the algebra $l^\infty (G_1)$ of all bounded functions on $G_1$, with the pointwise operations and $\sup$-norm.
\end{lem}

\medskip
{\it Proof}. It is obvious that $l^\infty (G_1)\subseteq M(A_1)$. Conversely, let $\rho_l\in L(A_1)$. Fix $s\in G_1$. There exists such a function $a\in A_1$ that $a(s)=1$. Denote $f(s)=\rho_l(a)(s)$. This element does not depend on the choice of $a$ -- if also $a'(s)=1$, then
$$
f(s)=f(s)a'(s)=\rho_l(a)(s)a'(s)=\rho_l(aa')(s)=\rho_l(a')(s)a(s)=\rho_l(a')(s).
$$
Therefore we may consider $s$ as a variable obtaining the function $f_l:\ G_1\to \K$.

For $b\in A_1$,
$$
\rho_l(b)(s)=\rho_l(b)(s)a(s)=\rho_l(a)(s)b(s)=f_l(s)b(s),
$$
so that $\rho_l$ is the operator of multiplication by $f_l$.

Note that $f_l\in l^\infty (G_1)$. Otherwise there would exist such a sequence $\{ g_n\}\subset G_1$ that $0\ne |f_l(g_n)|\to \infty$. Choose a function $b\in A_1$, equal to 0 everywhere outside this sequence and equal to $\dfrac1{f_l(g_n)}$ on it. Then the function $f_lb$ does not belong to $c_0(G_1)$, and we have come to a contradiction.

Similarly, every right multiplier is an operator of multiplication by a function $f_r\in l^\infty (G_1)$. The consistency condition (2.5) means that $f_rab=f_lab$ for any $a,b\in A_1$, so that $f_r=f_l$, thus $M(A_1)\subseteq l^\infty (G_1)$. $\qquad \blacksquare$

In order to define a comultiplication, we need a description of $M(A_1\tp A_1)$.

\medskip
\begin{lem}
There are the isomorphisms $A_1\tp A_1=c_0(G_1\times G_1)$, $M(A_1\tp A_1)=l^\infty (G_1\times G_1)$.
\end{lem}

\medskip
{\it Proof.} For any $\K$-Banach space $X$, the space $A_1\tp X$ is isomorphic to the space of vector-valued sequences $c_0(G_1,X)$ (see \cite{Se}). For $X=A_1$, this space consists of functions of two variables $f(s,t)$, $s,t\in G_1$, such that: 1) for any $\epsilon >0$, there exists a finite set $H_\epsilon \subset G_1$, for which $|f(s,t)|<\epsilon$ for all $s\in G_1$, $t\in G_1\setminus H_\epsilon$; 2) for every $t\in G_1$, there exists such a finite set $S_{t,\epsilon}\subset G_1$ that $|f(s,t)|<\epsilon$ for $s\in G_1\setminus S_{t,\epsilon}$, $t\in G_1$.

In particular, for $t\in H_\epsilon$, since the set $H_\epsilon$ is finite, one can choose $S_{t,\epsilon}=S_\epsilon$ independent of $t$.

Therefore, for $t\notin H_\epsilon$ we have $|f(s,t)|<\epsilon$ for all $s\in G_1$, while for $t\in H_\epsilon$ we have $|f(s,t)|<\epsilon$ for all $s\notin S_\epsilon$. In other words, $|f(s,t)|<\epsilon$ for $(s,t)\notin S_\epsilon\times H_\epsilon$. This means that $A_1\tp A_1=c_0(G_1\times G_1)$. By Lemma 5.1, $M(A_1\tp A_1)=l^\infty (G_1\times G_1)$. $\qquad \blacksquare$

\medskip
Define a comultiplication $\Delta:\ c_0(G_1)\to l^\infty (G_1\times G_1)$ setting
\begin{equation}
(\Delta f)(s,t)=f(st),\quad f\in c_0(G_1), s,t\in G_1.
\end{equation}
Then, for any $a,b\in A_1$, $s,t\in G_1$,
\begin{equation}
T_1(a\otimes b)(s,t)=T_3(a\otimes b)(s,t)=a(st)b(t);
\end{equation}
\begin{equation}
T_2(a\otimes b)(s,t)=T_4(a\otimes b)(s,t)=a(s)b(st).
\end{equation}
It is obvious that the right-hand sides in (5.2) and (5.3) belong to $c_0(G_1\times G_1)$.

Moreover, the above mappings are isometric. To prove the isometry property of an operator, it is sufficient to check that the operator transforms an orthonormal basis into an orthonormal basis. Such a basis of $A_1=c_0(G_1)$ is formed by functions $\delta_\xi (t)=\begin{cases}
1, & \text{if $t=\xi$},\\
0, & \text{if $t\ne \xi$}.
\end{cases}$.
An orthonormal basis in $A_1\tp A_1=c_0(G_1\times G_1)$ is given by the system of functions $(s,t)\mapsto \delta_\xi (s)\delta_\eta (t)$, $(\xi ,\eta )\in G_1\times G_1$, that is the system $\{ \delta_\xi \otimes \delta_\eta\}$.

We have
$$
T_1(\delta_\xi \otimes \delta_\eta )(s,t)=\delta_\xi (st) \delta_\eta (t)=\delta_{\xi \eta^{-1}}(s)\delta_\eta (t),
$$
so that the image under $T_1$ of the above orthonormal basis is a rearrangement of the latter. Similarly we check the isometry property of $T_2$. This reasoning proves also the surjectivity of these mappings.

As in \cite{VD98}, it is easy to calculate that
\begin{equation}
\varepsilon (f)=f(e),\quad (S(f))(t)=f(t^{-1}), \quad f\in A_1,t\in G_1,
\end{equation}
where $e$ is the identity element in $G_1$. Note that $a\delta_e=a(e)\delta_e=\varepsilon (a)\delta_e$ for all $a\in A_1$, that is $A_1$ is a multiplier Banach-Hopf algebra of discrete type, and its dual algebra is unital.

The extension of the above homomorphisms to multipliers is given explicitly by the same formulas written for $f\in l^\infty (G_1)$. However the general results are applicable too. We have only to construct a bounded approximate identity $\{ e_\lambda\}_{\lambda \in \Lambda}$.

Let $\Lambda$ be the set of all finite subsets of $G_1$ ordered by inclusion. For any $\lambda \in \Lambda$, define a function $e_\lambda \in A_1$ setting
$$
e_\lambda (s)=\begin{cases}
1, & \text{if $s\in \lambda$},\\
0, & \text{if $s\notin \lambda$}.
\end{cases}
$$
For any $f\in A_1$,
$$
(fe_\lambda) (s)=\begin{cases}
f(s), & \text{if $s\in \lambda$},\\
0, & \text{if $s\notin \lambda$},
\end{cases}
$$
so that
\begin{equation}
(fe_\lambda) (s)-f(s)=\begin{cases}
0, & \text{if $s\in \lambda$},\\
-f(s), & \text{if $s\notin \lambda$},
\end{cases}
\end{equation}

For every $\epsilon >0$, there exists such $\lambda_0\in \Lambda$ that $|f(s)|<\epsilon $ for $s\notin \lambda$, if $\lambda \succeq \lambda_0$. By (5.5), $|(fe_\lambda) (s)-f(s)|<\epsilon$ for all $s\in G_1$, if $\lambda \succeq \lambda_0$, so that $fe_\lambda \to f$.

Following \cite{VD98} once more, we find that the left- and simultaneously right-invariant functional on $A_1$ is given by the equality
\begin{equation}
\varphi (f)=\psi (f)=\sum\limits_{s\in G_1}f(s),\quad f\in A_1.
\end{equation}
Note however that the convergence in (5.6) is a purely non-Archimedean phenomenon; see Theorem 2.5.1 in \cite{PGS} regarding the summation of possibly uncountable sequences.

In order to check for this case the equality (4.1), note first that
$$
\sup\limits_{x\ne 0}\frac{|\varphi (xa)|}{\|x\|}\le \|a\|
$$
by the ultrametric inequality. On the other hand, $\|a\|=|a(s_0)|$ for some $s_0\in G_1$. If $x=\delta_{s_0}$, then $\varphi (xa)=a(s_0)$ and
$$
\sup\limits_{x\ne 0}\frac{|\varphi (xa)|}{\|x\|}\ge \frac{|\varphi (\delta_{s_0}a)|}{\|\delta_{s_0}\|}=|a(s_0)|=\|a\|,
$$
as desired.

Let us consider the dual object $\A_1$. Let $\omega_1=\varphi (\cdot a_1),\omega_2=\varphi (\cdot a_2)\in \A_1$. For any $b\in A_1$,
$$
(\omega_1\omega_2)(b)=(\omega_1\otimes \omega_2)\Delta (b)=\sum\limits_{s\in G_1}a_1(s)\sum\limits_{t\in G_1}a_2(t)b(st),
$$
so that
\begin{equation}
(\omega_1\omega_2)(b)=\sum\limits_{\tau\in G_1}\left[ \sum\limits_{s\in G_1}a_1(s)a_2(s^{-1}\tau )\right] b(\tau)
\end{equation}
and
\begin{equation}
\Hat\varphi (\omega_1\omega_2)=\sum\limits_{s\in G_1}a_1(s)a_2(s^{-1}).
\end{equation}

The convolution structure in (5.7) shows that $\varphi (\cdot \delta_e)$ is the unit in $\A_1$. The equality (5.8) implies the norm reproducing property for $\A_1$. The proof is similar to the above proof for the case of $A_1$.

By the definition of comultiplication in the dual object, we find that for every $x,y\in A_1$,
\begin{equation}
(\D (\omega_1)(1\otimes \omega_2))(x\otimes y)=\sum\limits_{s,t\in G_1}a_1(s)a_2(s^{-1}t)x(s)y(t),
\end{equation}
so that in the correspondence of $\omega_j\in \A_1$ with $a_j\in A_1$, $\D (\omega_1)(1\otimes \omega_2)$ corresponds to $a_1(s)a_2(s^{-1}t)$. A similar formula can be written for right multipliers. We will return to these formulas in Section 5.3 below.

\bigskip
{\bf 5.2.} {\it Zero-dimensional groups}. Let $G_2$ be a zero-dimensional Hausdorff locally compact topological group. We assume that $G_2$ carries a $\K$-valued left-invariant measure $\mu_l$. By \cite{MS,vR}, this happens if either the residue field of $\K$ has characteristic zero, or this characteristic equals $p\ne 0$ and $G_2$ has a $p$-free compact open subgroup $O$, that is no open subgroup in $O$ has an index divisible by $p$. The well-known example: $\K =\mathbb Q_p$, $G_2=\mathbb Q_l$ with $l\ne p$. Under the same conditions, a $\K$-valued right-invariant measure $\mu_r$ exists too.

The commutative Banach algebra $A_2=C_0(G_2)$, with the $\sup$-norm and pointwise operations, consists of continuous functions $f:\ G_2\to \K$, such that for any $\epsilon >0$, the set $\{ s\in G_2:\ |f(s)|\ge \epsilon\}$ is compact.

\medskip
\begin{lem}
The multiplier algebra $M(A_2)$ is isomorphic to the algebra $C_b(G_2)$ of all bounded continuous functions from $G_2$ to $\K$ with the $\sup$-norm and pointwise operations.
\end{lem}

\medskip
{\it Proof}. As in the proof of Lemma 5.1, for each $\rho_l\in L(A_2)$ we find a function $f_l:\ G_2\to \K$, such that $(\rho_l(b))(s)=f_l(s)b(s)$, for all $b\in A_2$, $s\in G_2$. For a fixed $s$, $f_l(s)=\rho_l(a)(s)$ where $a\in A_2$ is chosen in such a way that $a(s)=1$. The existence of such a function $a\in A_2$ is proved in \cite{PGS}, Theorem 2.5.32.

The function $f_l$ is continuous as a ratio of two continuous functions. Let us prove its boundedness (we follow the method from \cite{Wa}). For any point $s\in G_2$, we set
$$
C_s=\sup\limits_{0\ne \gamma \in A_2}\frac{|\gamma (s)|}{\|\gamma \|}=\sup\limits_{\gamma \in A_2,\|\gamma\|=1}|\gamma (s)|
$$
(the above expressions are equal, since in this case the $\sup$-norm takes the same values as the absolute value $|\cdot |$). Then $0\le C_s\le 1$ and $|\gamma (s)|\le C_s\|\gamma \|$. In particular,
$$
|f_l(s)b(s)|\le C_s\|f_lb\|=C_s\|\rho_l(b)\|\le C_s\|\rho_l\|\cdot \|b\|
$$
for any $b\in A_2$, in particular, for any $b$ with $\|b\|=1$. Therefore
$$
|f_l(s)|\le \inf\limits_{\|b\|=1}\frac{C_s\|\rho_l\|}{|b(s)|}=\frac{C_s\|\rho_l\|}{\sup\limits_{\|b\|=1}|b(s)|}=\|\rho_l\|,
$$
so that $f_l$ is bounded, $\|f_l\|\le \|\rho_l\|.\qquad \blacksquare$

\medskip
It is known \cite{KB} that $A_2\tp A_2=C_0(G_2\times G_2)$. Therefore $M(A_2\tp A_2)=C_b(G_2\times G_2)$.

The comultiplication, counit and antipode are defined by the same formulas (5.1) and (5.4) as in the case of discrete groups. The isometry of the mappings (5.2) and (5.3) are proved in this case as follows.

It is known (\cite{PGS}, Theorem 2.5.34) that the Banach space $C_0(G_2)$ has an orthonormal basis $\{ \chi_\nu \}$ consisting of characteristic functions of open compact sets. Then $\{ \chi_\nu (s)\chi_\varkappa (t)\}$ is an orthogonal basis in $C_0(G_2\times G_2)$; see \cite{PGS}, Corollary 10.2.10 and Theorem 10.3.16. It is obvious that the characteristic functions of two sets are orthogonal if and only if neither of the sets contains the other. It is straightforward to check that the mappings $T_1=T_3$ and $T_2=T_4$ given by (5.2) and (5.3) maintain this property. Thus, they transform the above orthonormal basis into an orthonormal system of functions, which implies the isometry property.

To prove the surjectivity of $T_1$ ($T_2$ is considered in a similar way), note that any function of $x$ and $y$ from $C_0(G_2\times G_2)$ can be approximated uniformly as $\lim\sum a_i(x)b_i(y)$. In particular, we may write
$$
u(xy^{-1},y)=\lim\sum a_i(x)b_i(y)
$$
(note that the mapping $(s,t)\mapsto (st,t)$ is continuous, hence this function belongs to $C_0(G_2\times G_2)$), and then set $x=st,y=t$, $t\in G_2$, which results in the uniform limit
$$
u(s,t)=\lim\sum T_1(a_i\otimes b_i)(s,t).
$$

The invariant functionals on $A_2$ are defined as follows:
$$
\varphi (f)=\int\limits_{G_2} f(s)\mu_l(ds);\quad \psi (f)=\int\limits_{G_2} f(s)\mu_r(ds).
$$
It is known \cite{MS} that the $\K$-valued Haar measure can be normalized in such a way that $\|\varphi\| \le 1$ and the identity (4.1) is valid. The proofs of some properties of non-Archimedean integrals in \cite{MS} are given for integrands with compact supports, but are easily extended using standard approximation arguments.

In order to use the above general results, we need the existence of bounded approximate identities.

\medskip
\begin{lem}
The Banach algebras $A_2$ and $\A_2$ possess bounded approximate identities.
\end{lem}

\medskip
{\it Proof}. Note first of all that $G_2$ is paracompact (see \cite{HR}, Theorem 8.13), so that we can use Lemma 28.10 from \cite{DW} stating that there exists a family $H_\lambda$, $\lambda \in \Lambda$ ($\Lambda$ is a directed set), of compact open subsets of $G_2$, such that $H_{\lambda'}\subset H_{\lambda''}$ for $\lambda' \prec \lambda''$ and $G_2=\bigcup\limits_{\lambda \in \Lambda}H_\lambda$. Let $e_\lambda$ be the $\K$-valued characteristic function of $H_\lambda$. This function belongs to $A_2$, and for any $u\in A_2$,
$$
(ue_\lambda -u)(x)=\begin{cases}
0, & \text{for $x\in H_\lambda$},\\
-u(x), & \text{elsewhere}.
\end{cases}
$$

For any $\epsilon >0$, there exists such a compact set $F_\epsilon$ that $|u(x)|<\epsilon$, if $x\notin F_\epsilon$. The family $H_\lambda$ forms an open covering of $F_\epsilon$; there exists its finite subcovering, and since the family $H_\lambda$ is increasing, we find such $\lambda_\epsilon\in \Lambda$ that $F_\epsilon \subset H_\lambda$ for all $\lambda \succeq \lambda_\epsilon$. Therefore $|u(x)|<\epsilon$, if $x\notin H_\lambda$, so that $\|ue_\lambda -u\|<\epsilon$, if $\lambda \succeq \lambda_\epsilon$. Thus, $\{ e_\lambda \}$ is a bounded (by one) approximate identity in $A_2$.

The dual object $\A_2$ consists of functionals
$$
\omega_a (x)=\int\limits_{G_2}a(s)x(s)\mu_l(ds),\quad a\in A_2,
$$
with the multiplication
\begin{multline*}
\left( \omega_{a'}\omega_{a''}\right) (x)=\left( \omega_{a'}\otimes \omega_{a''}\right) \Delta (x)=\int\limits_{G_2}a'(s)\mu_l(ds)\int\limits_{G_2}a''(t)x(st)\mu_l(dt)\\
=\int\limits_{G_2}x(\theta )(a'\star a'')(\theta )\mu_l(d\theta)
\end{multline*}
where the convolution $a'\star a''$ has the form
$$
(a'\star a'')(\theta )=\int\limits_{G_2}a'(s)a''(s^{-1}\theta)\mu_l(ds)
$$
(for an analog of Fubini's theorem for $\K$-valued measures see \cite{MS}).

With this multiplication, $\A_2$ is isomorphic to the group Banach algebra $L(G_2)$ studied in \cite{vRS} where, in particular, a bounded approximate identity for this algebra is constructed. Note that the right multiplier algebra of $L(G_2)$ is isomorphic to the algebra of improper measures on $G_2$; see Exercise 8.B.v in \cite{vR}.

As we know, $\Hat \varphi (\omega_a)=\varepsilon (a)=a(e)$, so that
$$
\Hat \varphi (\omega_{a'}\omega_{a''})=\int\limits_{G_2}a'(s)a''(s^{-1})\mu_l(ds).
$$
The norm reproducing property for $\A_2$ follows from the result from \cite{MS} already used above.

\bigskip
{\bf 5.3.} {\it Algebras generated by regular representations}. Let $G_3$ be a discrete group. On the $\K$-Banach space $c_0(G_3)$, we consider the right regular representation of $G_3$: $(R_af)(s)=f(sa)$, $s,a\in G_3$. Denote by $\RR$ the set of all the operators $R_a,a\in G_3$. The closure of the linear span of $\RR$ in the strong operator topology coincides with its closure in the uniform operator topology and is isomorphic to the Banach algebra $A_3=A(G_3)$ equal, as a Banach space, to $c_0(G_3)$ but, in contrast to $A_1$, endowed with the product
\begin{equation}
(u*v)(d)=\sum\limits_{l\in G_3}u(l)v(l^{-1}d),\quad d\in G_3;
\end{equation}
for the details see \cite{K14}.

As we know (Section 5.1), $\{ \delta_s\}_{s\in G_3}$ is an orthonormal basis in $A_3$. The element $\delta_e$, where $e$ is the unit in $G_3$, is the unit in $A_3$, so that $A_3$ is equal to its multiplier algebra.

Define a comultiplication in $A_3$ setting $\Delta (\delta_s)=\delta_s \otimes \delta_s$ and extending this by linearity and continuity. In other words, writing an arbitrary $u\in A_3$ as
$$
\Delta (u)=\sum\limits_{s\in G_3}u(s)\delta_s
$$
we have
$$
u=\sum\limits_{s\in G_3}u(s)\delta_s \otimes \delta_s.
$$

By Lemma 5.2, $A_3\tp A_3$ can be identified with $c_0(G_3\times G_3)$ as a Banach space; this identification extends to the identity of Banach algebras, if we endow $c_0(G_3\times G_3)$ with the coordinate-wise convolution product. In this setting, for any $u\in A_3$,
$$
(\Delta (u))(s,t)=\begin{cases}
u(s), & \text{if $s=t$};\\
0, & \text{if $s\ne t$},
\end{cases}
\quad s,t\in G_3,
$$
so that
\begin{equation}
T_1(u\otimes v)(s,t)=(\Delta (u)(1\otimes v))(s,t)=\sum\limits_{l\in G_3}(\Delta (u))(s,l)v(l^{-1}t)=u(s)v(s^{-1}t).
\end{equation}

Obviously, $T_1(u\otimes v)$ belongs to $A_3\tp A_3$. Since any function $F(s,st)$, $F\in A_3\tp A_3$, can be approximated uniformly by the expressions $\sum c_ju_j(s)v_j(t)$, $c_j\in \K$, $u_j,v_j\in A_3$, for $t=s^{-1}\tau$ this gives the uniform approximation of any function $F(s,\tau )$, thus the surjectivity of $T_1$.

More generally, we can write
$$
T_1(F)(s,t)=F(s,s^{-1}t)
$$
for any $F\in A_3\tp A_3$. This implies, in particular, that $T_1$ is an isometry. Similarly, $T_2,T_3$, and $T_4$ can be considered.

It is easy to find the counit and antipode:
$$
\varepsilon (u)=\sum\limits_{s\in G_3}u(s),\quad S(u)(t)=u(t^{-1}),
$$
and the invariant functionals:
$$
\varphi (u)=\psi (u)=u(e).
$$
The equality (4.1) is checked just as its counterpart for the algebra $A_1$.

Comparing (5.10) with (5.7), and (5.11) with (5.9), we see that $A_3$ is isomorphic to $\A_1$.


\medskip
\begin{thebibliography}{999}
\bibitem{BGR}
S. Bosch, U. G\"untzer, and R. Remmert, {\it Non-Archimedean
Analysis}, Springer, Berlin, 1984.
\bibitem{Da}
M. Daws, Multipliers, self-indiced and dual Banach algebras, {\it Dissertationes Math.} {\bf 470} (2010), 62 pp.
\bibitem{D74}
B. Diarra, La dualit\'e de Tannaka dans les corps valu\'es ultram\'etriques complets, {\it C. R. Acad. Sci. Paris, S\'er. A} {\bf 279}, No. 26 (1974), 907--909.
\bibitem{D79}
B. Diarra, Sur quelques repr\'esentations p-adiques de $\mathbb Z_p$, {\it Indag. Math.} {\bf 41} (1979), 481--493.
\bibitem{D96}
B. Diarra, Alg\`ebres de Hopf et fonctions presque p\'eriodiques ultram\'etriques, {\it Rivista di Mat. Pura ed Appl.} {\bf 17} (1996), 113--132.
\bibitem{DW}
R. S. Doran and J. Wichmann, {\it Approximate Identities and Factorization in Banach Modules}, Lect. Notes Math. 768 (1979), 305 pp.
\bibitem{HR}
E. Hewitt and K. A. Ross, {\it Abstract Harmonic Analysis}, Vol. 1, Springer, Berlin, 1963.
\bibitem{Jo}
B. E. Johnson, An introduction to the theory of centralizers, {\it Proc. London Math. Soc.} {\bf 14} (1964), 299--320.
\bibitem{KB}
A. K. Katsaras and A. Beloyiannis, Tensor products of non-Archimedean  weighted spaces of continuous functions, {\it Georgian Math. J.} {\bf 6}, No. 1 (1999), 33--44.
\bibitem{Ki}
J. Kisy\'nski, On Cohen's proof of the Factorization Theorem, {\it Ann. Polon. Math.} {\bf 75}, No. 2 (2000), 177--192.
\bibitem{K13}
A. N. Kochubei,  On some classes of non-Archimedean operator algebras, {\it Contemporary Math.} {\bf 596} (2013), 133--148.
\bibitem{K14}
A. N. Kochubei, Non-Archimedean group algebras with Baer reductions, {\it Algebras and Represent. Theory} {\bf 17} (2014), 1861--1867.
\bibitem{KV00}
J. Kustermans and S. Vaes, Locally compact quantum groups, {\it Ann. Sci. \'Ecole Norm. Sup.} {\bf 33} (2000), 837--934.
\bibitem{KV03}
J. Kustermans and S. Vaes, Locally compact quantum groups in the von Neumann algebraic setting, {\it Math. Scand.} {\bf 92} (2003), 68--92.
\bibitem{MS}
A. F. Monna and T. A. Springer, Int\'egration non-archim\'edienne, {\it Indag. Math.} {\bf 25} (1963), 634--653.
\bibitem{PGS}
C. Perez-Garcia and W. H. Schikhof, {\it Locally Convex Spaces over
Non-Archimedean Valued Fields}, Cambridge University Press, Cambridge,
2010.
\bibitem{vR}
A. C. M. van Rooij, {\it Non-Archimedean Functional Analysis}, Marcel Dekker, New York, 1978.
\bibitem{vRS}
A. C. M. van Rooij and W. H. Schikhof, Group representations in non-Archimedean Banach spaces, {\it Bull. Soc. Math. France. Supplement}, M\'emoire 39-40 (1974), 329--340.
\bibitem{Sch1}
W. H. Schikhof, {\it Non-Archimedean Harmonic Analysis}. Thesis, Nijmegen, 1967.
\bibitem{Sch2}
W. H. Schikhof, Non-Archimedean representations of compact groups, {\it Compositio Math.} {\bf 23} (1971), 215--232.
\bibitem{Se}
J.-P. Serre, Endomorphismes compl\`etement continus des espaces de
Banach $p$-adiques, {\it Publ. Math. IHES} {\bf 12} (1962),
69--85.
\bibitem{So}
Y. Soibelman, Quantum p-adic spaces and quantum p-adic groups. In: {\it ``Geometry and Dynamics of Groups and Spaces''}, Progr. Math. 265, Birkh\"auser, Basel, 2008, pp. 697--719.
\bibitem{Sun}
Menghong Sun, {\it Multiplier Hopf Algebras and Duality}. Thesis. University of Windsor, Ontario, Canada, 2011. http://scholar.uwindsor.ca/etd/352 .
\bibitem{Ti}
Th. Timmermann, {\it An Invitation to Quantum Groups and Duality}. European Mathematical Society, Z\"urich, 2008.
\bibitem{VVD}
S. Vaes and A. Van Daele, Hopf $C^*$-algebras, {\it Proc. London Math. Soc.} {\bf 82} (2001), 337--384.
\bibitem{VD94}
A. Van Daele, Multiplier Hopf algebras, {\it Trans. Amer. Math. Soc.} {\bf 342} (1994), 917--932.
\bibitem{VD98}
A. Van Daele, An algebraic framework for group duality, {\it Adv. Math.} {\bf 140} (1998), 323--366.
\bibitem{VD14}
A. Van Daele, Locally compact quantum groups. A von Neumann algebra approach, {\it SIGMA} {\bf 10} (2014), article 082.
\bibitem{VDW}
A. Van Daele and S. Wang, Weak multiplier Hopf algebras. Preliminaries, motivation and basic examples. {\it Banach Center Publications} {\bf 98} (2012), 367--415.
\bibitem{Wa}
J.-K. Wang, Multipliers on commutative Banach algebras, {\it Pacif. J. Math.} {\bf 11} (1961), 1131--1149.
\end{thebibliography}
\end{document}